
\documentstyle{amsppt}
\baselineskip18pt
\magnification=\magstep1
\pagewidth{30pc}
\pageheight{45pc}

\hyphenation{co-deter-min-ant co-deter-min-ants pa-ra-met-rised
pre-print pro-pa-gat-ing pro-pa-gate
fel-low-ship Cox-et-er dis-trib-ut-ive}
\def\leaderfill{\leaders\hbox to 1em{\hss.\hss}\hfill}

\def\Be{{\Cal B}}
\def\C{{\Cal C}}

\def\L{{\Cal L}}

\def\id{\text {\rm id}}

\def\idest{i.e.,\ }

\def\a{{\alpha}}
\def\be{{\beta}}
\def\g{{\gamma}}

\def\d{{\delta}}

\def\rreal{\Delta^{re}}
\def\preal{\Delta^{re}_+}
\def\e{{\varepsilon}}
\def\z{{\zeta}}

\def\l{{\lambda}}

\def\s{{\sigma}}

\def\cha{{\chi}}
\def\w{{\omega}}

\def\bs{{\bold s}}

\def\pmlp{\bold{2^{l+1}}}
\def\pml{\bold{2^l}}
\def\pmlm{\bold{2^{l-1}}}
\def\pmc{\bold{2}_+^{\zed_{2l}}}

\def\b0{\text{\bf 0}}

\def\ra{{\ \longrightarrow \ }}

\def\lan{{\langle}}
\def\ran{{\rangle}}

\def\lrt#1#2{\left\langle {#1}, {#2} \right\rangle}

\def\lpr{\lan \phi \ran}

\def\fg{{\frak g}}

\def\complex{{\Bbb C}}
\def\zed{{\Bbb Z}}

\def\boxit#1{\vbox{\hrule\hbox{\vrule \kern3pt
\vbox{\kern3pt\hbox{#1}\kern3pt}\kern3pt\vrule}\hrule}}
\def\rabbit{\vbox{\hbox{\kern0pt
\vbox{\kern0pt{\hbox{---}}\kern3.5pt}}}}

\def\tableau#1{
        \hbox {
                \hskip -10pt plus0pt minus0pt
                \raise\baselineskip\hbox{
                \offinterlineskip
                \hbox{#1}}
                \hskip0.25em
        }
}

\def\tabCol#1{
\hbox{\vtop{\hrule
\halign{\strut\vrule\hskip0.5em##\hskip0.5em\hfill\vrule\cr\lower0pt
\hbox\bgroup$#1$\egroup \cr}
\hrule
} } \hskip -10.5pt plus0pt minus0pt}

\def\CR{
        $\egroup\cr
        \noalign{\hrule}
        \lower0pt\hbox\bgroup$
}



\def\blank#1#2{
\hbox to #1{\hfill \vbox to #2{\vfill}}
}


\def\strut{\vrule height10pt depth5pt width0pt}

\def\secaa{1}
\def\secab{2}
\def\secba{3}
\def\secbb{4}
\def\sectc{5}
\def\sectd{6}
\def\secea{7}
\def\seceb{8}
\def\sectf{9}
\def\esla{a}
\def\eslb{b}
\def\eslc{c}
\def\apage{\page}

\topmatter
\title Full heaps and representations of affine Weyl groups
\endtitle

\author R.M. Green \endauthor
\affil Department of Mathematics \\ University of Colorado \\
Campus Box 395 \\ Boulder, CO  80309-0395 \\ USA \\ {\it  E-mail:}
rmg\@euclid.colorado.edu \\
\newline
\endaffil

\abstract
We use the author's combinatorial theory of full heaps to categorify the
action of a large class of Weyl groups on their root systems, and thus to
give an elementary and uniform construction of a family of faithful 
permutation representations of Weyl groups.  Examples
include the standard representations of affine Weyl groups as permutations
of $\zed$ and geometrical examples such as the realization of the Weyl
group of type $E_6$ as permutations of $27$ lines on a cubic surface; in
the latter case, we also show how to recover the 
incidence relations between the
lines from the structure of the heap.  Another class of examples involves the
action of certain Weyl groups on sets of pairs $(t, f)$, where
$t \in \zed$ and $f$ is a function from a suitably chosen set to the
two-element set $\{+, -\}$.
Each of the permutation representations corresponds to a module for
a Kac--Moody algebra, and gives an explicit basis for it.
\endabstract

\toc
\head{} Introduction \apage{2}\endhead
\head{1.} Heaps over Dynkin diagrams \apage{6}\endhead
\head{2.} Full heaps \apage{9}\endhead
\head{3.} The Weyl group and skew ideals \apage{12}\endhead
\head{4.} Main results \apage{18}\endhead
\head{5.} Permutations of $\zed$ \apage{22}\endhead
\head{6.} Geometrical examples \apage{31}\endhead
\head{7.} The binary path representation in type affine $B$ \apage{36}\endhead
\head{8.} Other binary path representations \apage{40}\endhead
\head{9.} Some related constructions \apage{44}\endhead
\head{} Acknowledgement \apage{48}\endhead
\head{} References \apage{48}\endhead
\endtoc
\endtopmatter
\centerline{\bf Preliminary version, draft 3}

\head Introduction \endhead

In \cite{{\bf 7}}, we introduced the notion of full heaps, which are remarkable
locally finite labelled posets that are closely related to to Kac--Moody 
algebras, distributive lattices, crystal bases and Weyl groups.  The main 
result of \cite{{\bf 7}} is a construction of almost all affine Kac--Moody 
algebras modulo their one-dimensional centres in terms of raising and 
lowering operators on the space spanned by 
the so-called proper ideals of a suitable full heap.  We mentioned briefly
in \cite{{\bf 7}, \S8} that there is a natural action of the Weyl group $W$ on
the set $\Be$ of proper ideals of a full heap $E$, and the purpose of this 
paper is to understand this action.

We concentrate on the case where the full heap corresponds
to an affine Kac--Moody algebra.
Our strategy here is to use the distributive lattice $\Be$ to categorify
the root system associated to $W$, in a way that is compatible with the
action of $W$; this is based on the categorification of positive roots 
implicit in \cite{{\bf 7}}.
More precisely, we define a ``decategorification'' map $\cha$ from 
$\Be \times \Be$
taking values in the root lattice of $W$, and we concentrate on those
pairs $(F, F') \in \Be \times \Be$ whose character is a real root; in
such a situation we say that $F$ is ``skew'' to $F'$.
We prove in Theorem \secbb.1 that the diagonal action of $W$ on
$\Be \times \Be$ induces, via $\cha$, the usual action of $W$ on the root
system.  The standard fact that $W$ acts faithfully on the (real) roots then
implies that the action of $W$ on $\Be$ is faithful.

The permutation representations of Weyl groups arising from the action of
$W$ on $\Be$ are very interesting, and we give explicit descriptions of
some of the most important ones in this paper.  Loosely speaking, these
representations seem to be of three types, depending on the structure of the
full heap $E$ that gives rise to them, although the definitions themselves
are completely uniform.

The first type of permutation representation arises from full heaps
all of whose antichains are short (\idest usually of size $1$, and otherwise
of size $2$), and the permutation
representations arising in this way include the standard realizations of
the affine Weyl groups of types $A$, $B$, $C$ and $D$ as permutations of
the integers.  In this case, most pairs $(F, F')$ of elements
of $\Be$ will be skew, and it is this close relationship between the lattice
$\Be$ and the root system that makes these representations so useful for
understanding the Coxeter group structure of $W$.  A comprehensive guide
to these representations (and their applications) may be found in 
\cite{{\bf 3}, \S8}, and we give a brief history of them at the beginning 
of \S\sectc.  These representations also arise 
(together with a representation of the affine Weyl group of type $G_2$) 
in recent work of Cellini,
M\"oseneder Frajria and Papi \cite{{\bf 4}, \S4}, as a by-product of their work 
on a combinatorial interpretation of Kostant's formula for powers of the 
Euler product.

The second type of permutation representation arises from full heaps all of
whose antichains are long, \idest of size approaching $n/2$, where $n$ is
the number of elements in the Dynkin diagram.  In these cases, the heap $E$
is as far from being totally ordered as possible, so labelling the elements
of $\Be$ by integers is not convenient and depends heavily on choice.
However, we will show how to parametrize the elements of 
$\Be$ by certain pairs $(t, f)$, where $t$ is an integer
and $f$ is a function from some set (depending on $E$) to the set
$\{+, -\}$ with two elements.  For permutation representations of this
type, most pairs $(F, F')$ of
elements of $\Be$ will not be skew.  These examples, which are described in
\S\secea\  and \S\seceb, are new from the Weyl group point of view
to the best of our knowledge, although they are reminiscent of the wreath
product constructions of the finite Weyl groups of types $B$ and $D$
\cite{{\bf 9}, \S2.10}. 

The third type of permutation representation includes all other cases, meaning
that the antichains of the heap are of intermediate length.  
The relation of skewness on $\Be$ is the most interesting
in these cases, but some ingenuity may be required to obtain an appropriate 
parametrization.  We look at one such example
in detail in \S\sectd, namely the case of type $E_6$, which we obtain from type
$E_6^{(1)}$ by restriction.  The full heap construction realizes the finite
Weyl group as a permutation group on $27$ objects.  These objects can
be identified in a natural way with the $27$ lines on a cubic surface, and
remarkably, the combinatorial notion of skewness coincides with the
geometric notion of skewness on the $27$ lines.  The approach here also makes
it obvious how to lift the action of $E_6$ on the $27$ lines to the affine
Weyl group of type $E_6^{(1)}$ in an explicit way (Proposition \sectd.1).
It is already known that the $27$ lines are in correspondence with the
weights of a minuscule representation of the Lie algebra $E_6$ and that 
two lines are incident if and only if the corresponding weights are not
orthogonal with respect to a certain inner product (see \cite{{\bf 13}, \S1, \S3}).
However, the full heaps approach is more elementary 
in that one need only use the theory of Coxeter groups and their root 
systems, and the representation theory of Lie algebras is not required
to describe the construction.  
There ought to be a somewhat similar geometric construction for type $E_7$, 
but details of the correspondence have yet to be worked out.  There are
other examples of this third type of representation, including some
associated to the Coxeter system of type $A_l^{(1)}$, but since they are not
yet well understood, we will not consider them here.

The injective homomorphisms from affine Weyl groups to permutations of $\zed$
used by Cellini et al \cite{{\bf 4}} and by Eriksson \cite{{\bf 6}}
are defined in terms of the action of the affine Weyl group on a carefully 
chosen vector.  An advantage of our combinatorial point of view using 
full heaps over these two approaches is that it is extremely elementary: 
the only Lie theory needed for 
our construction is the definition of a Dynkin diagram, or generalized 
Cartan matrix, and we do not even need any linear algebra.
However, the proof that the map we define is indeed an 
injective homomorphism does use Lie theoretic concepts, such as 
the theory of Coxeter groups.

As we explained in \cite{{\bf 7}}, a full heap over the Dynkin diagram of an 
affine Kac--Moody algebra determines on the one hand a representation of
the Kac--Moody algebra and on the other hand a representation
of the corresponding affine Weyl group.  The permutation representations we
study in this paper are therefore closely related to the representation
theory of Lie algebras.  In particular, the parametrizations of the sets 
$\Be$ given in this paper have immediate applications to Lie algebras.  
For example, it follows from Proposition \secea.8 (respectively, Proposition
\seceb.2) that the finite Weyl group of type $B_n$ (respectively, $D_n$)
has a natural faithful action on the set of all strings of length $n$ 
from the alphabet $\{+, -\}$ (respectively, the set of all strings of length
$n$ from the alphabet $\{+, -\}$ that contain an even number of occurrences
of $-$).  From the Lie algebra point of view described in \cite{{\bf 7}}, these
sets of strings parametrize crystal bases of spin representations
of the corresponding simple Lie algebras over $\complex$, and the action 
of a Chevalley basis on them may be
explicitly calculated.  Note that we do not need Clifford algebras to do
this, and the heaps approach makes it obvious why the modules have 
dimensions $2^n$ and $2^{n-1}$, respectively.    
Although this result could also be achieved directly using
the theory of crystal bases, the full heap may be a much simpler object than
the corresponding crystal (see remarks \seceb.8 and \sectf.3).

\vfill\eject

\head \S\secaa. Heaps over Dynkin diagrams \endhead

We first review from \cite{{\bf 7}, \S1} some of the basic properties of
heaps over Dynkin diagrams.
The definitions relating to generalized Cartan matrices come
from \cite{{\bf 10}}, and the heap definitions are based on \cite{{\bf 18}}.

Let $A$ be an $n$ by $n$ matrix with integer entries.
We call $A$ a {\it generalized Cartan matrix} 
if it satisfies the conditions (a) $a_{ii} = 2$
for all $1 \leq i \leq n$, (b) $a_{ij} \leq 0$ for $i \ne j$ and (c)
$a_{ij} = 0 \Leftrightarrow a_{ji} = 0$.  In this paper, we will only
consider generalized Cartan matrices with entries in the set 
$\{2, 0, -1, -2\}$; such matrices are sometimes called {\it doubly laced}.  
If, furthermore, $A$ has no entries equal to $-2$, 
we will call $A$ {\it simply laced}.  

The Dynkin diagram $\Gamma = \Gamma(A)$ associated to a generalized Cartan
matrix is a directed graph, possibly with multiple edges, and vertices
indexed (for now) by the integers $1$ up to $n$.  If $i \ne j$ and
$|a_{ij}| \geq |a_{ji}|$, 
we connect the vertices corresponding to $i$ and $j$ by $|a_{ij}|$ lines;
this set of lines is equipped with an arrow 
pointing towards $i$ if $|a_{ij}|
> 1$.  For example, if $a_{ij} = a_{ji} = -2$, this will result in a double
edge between $i$ and $j$ equipped with an arrow pointing in each direction.
There are further rules if $a_{ij} a_{ji} > 4$, but we do not need
these for our purposes.

The Dynkin diagram (together with the 
enumeration of its vertices) and the generalized Cartan matrix determine 
each other, so we may write $A = A(\Gamma)$.  If $\Gamma$ is 
connected, we call $A$ {\it indecomposable}.

Let $\Gamma$ be a Dynkin diagram with vertex set $P$ and no 
multiple edges.  Let $C$ be the relation on $P$ such that $x \ C\ y$ if and
only if $x$ and $y$ are distinct unadjacent vertices in $\Gamma$, and let
$\C$ be the complementary relation.

\definition{Definition \secaa.1}
A {\it labelled heap} over $\Gamma$ is a triple $(E, \leq, \e)$ 
where $(E, \leq)$ is a locally finite partially ordered set (in other words,
a poset all of whose intervals are finite) with order relation denoted
by $\leq$ and where $\e$ is a map $\e : E \ra P$ satisfying the following two
axioms. 

\item{1.}{For every $\a, \be \in E$ such that $\e(\a) \ \C \ \e(\be)$, 
$\a$ and $\be$ are comparable in the order $\leq$.}

\item{2.}{The order relation $\leq$ is the transitive closure of the
relation $\leq_\C$ such that for all $\a, \be \in E$, $\a \ \leq_\C \ \be$ 
if and only if both $\a \leq \be$ and $\e(\a) \ \C \ \e(\be)$.}
\enddefinition

We call $\e(\a)$ the {\it label} of $\a$.  In the sequel, we will sometimes
appeal to the fact that the partial order is the reflexive, transitive 
closure of the covering relations, because of the local finiteness condition.

\definition{Definition \secaa.2}
Let $(E, \leq, \e)$ and $(E', \leq', \e')$ be two labelled 
heaps over $\Gamma$.  We say that $E$ and $E'$ are isomorphic (as labelled
posets) if there is a poset isomorphism $\phi : E \ra E'$ such that 
$\e = \e' \circ \phi$.

A {\it heap} over $\Gamma$ is an isomorphism class of labelled
heaps.  We denote the heap corresponding to the labelled heap 
$(E, \leq, \e)$ by $[E, \leq, \e]$.
\enddefinition

We will sometimes abuse language and speak of the underlying set of a heap,
when what is meant is the underlying set of one of its representatives.

\definition{Definition \secaa.3}
Let $(E, \leq, \e)$ be a labelled heap over 
$\Gamma$, and let $F$ a subset of $E$.  
Let $\e'$ be the restriction of $\e$ to $F$.  Let ${\Cal R}$ be the relation
defined on $F$ by $\a \ {\Cal R} \ \be$ if and only if $\a \leq \be$ and
$\e(\a) \ \C \ \e(\be)$.  Let $\leq'$ be the transitive closure of ${\Cal R}$.
Then $(F, \leq', \e')$ is a labelled heap over $\Gamma$.  The heap 
$[F, \leq', \e']$ is called a {\it subheap} of $[E, \leq, \e]$.  

If $E = (E, \leq, \e)$ is a labelled heap over $\Gamma$, then we define
the {\it dual labelled heap}, $E^*$ of $E$, to be the labelled heap 
$(E, \geq, \e)$.  (The notion of ``dual heap'' is defined analogously.)
There is a natural anti-isomorphism of labelled posets from $E$ to $E^*$,
which we will denote by $*$.  

Recall that if $(E, \leq)$
is a partially ordered set, a function $\rho : E \ra \zed$ is said to be a 
{\it rank function} for $(E, \leq)$ if whenever $a, b \in E$ are such that
$a < b$ is a covering relation, we have $\rho(b) = \rho(a) + 1$.  If a
rank function for $(E, \leq)$ exists, we say $(E, \leq)$ is {\it ranked},
and we say that the heap $(E, \leq, \e)$ is ranked to mean that $(E, \leq)$
is ranked as a partially ordered set.

If $F$ is convex as a subset of $E$ (in other words, if 
$\a \leq \be \leq \g$ with $\a, \g \in F$, then $\be \in F$) then we call 
$F$ a {\it convex subheap} of $E$.  If, whenever $\a \leq \be$ and $\be 
\in F$ we have $\a \in F$, then we call $F$ an {\it ideal} of $E$.
If $F$ is an ideal of $E$ with $\emptyset \subsetneq F \subsetneq E$ 
such that for each vertex $p$ of $\Gamma$ we have 
$\emptyset \subsetneq F \cap \e^{-1}(p) \subsetneq \e^{-1}(p)$, then we
call $F$ a {\it proper ideal} of $E$.  If $\g_1, \g_2, \ldots, \g_k$ are
elements of $E$, then we define the ideal $$
\lan \g_1, \g_2, \ldots, \g_k \ran = \{\a \in E : \a \leq \g_i \text{ for
some } 1 \leq i \leq k\}
.$$  If $F = \lan \g \ran$ for some $\g \in E$, then we call $F$ a 
{\it principal} ideal.

We call $E$ {\it periodic} if there
exists a nonidentity automorphism $\phi : E \ra E$ of labelled posets
such that $\phi(x) \geq x$ for all $x \in E$.  
By \cite{{\bf 7}, Remark 7.2}, $\phi$ restricts to an automorphism of 
the chains $\e^{-1}(p)$ for $p$ a vertex of $\Gamma$ of the form 
$\phi(E(p, x)) = E(p, x + t_p)$ for some
nonnegative integer $t_p$ depending on $p$ but not on the labelling chosen
for $E$.  If $\a \in R^+$ is such that $\a(p) = t_p$, we will say that $\phi$ 
is {\it periodic with period $\a$}.  If there is no automorphism $\phi'$ of $E$
with period $\a'$ such that $\a = n \a'$ with $n > 1$, then we also say that
$E$ is periodic with period $\a$ and fundamental automorphism $\phi$.
\enddefinition

We will often implicitly use the fact that a subheap is determined by its
set of vertices and the heap it comes from.  Note that in a periodic heap
$E$, the automorphism $\phi$ induces an inclusion-preserving permutation
of the proper ideals of $E$.
 
\definition{Definition \secaa.4}
Let $(E, \leq, \e)$ be a locally finite labelled heap over $\Gamma$.
We say that $(E, \leq, \e)$ and $[E, \leq, \e]$ are {\it fibred}
if
\item{(a)}{for each vertex $p$ in $\Gamma$, the subheap $\e^{-1}(p)$ is
unbounded above and unbounded below,}
\item{(b)}{for every pair $p, p'$ of adjacent vertices in $\Gamma$
and every element $\a \in E$ with $\e(\a) = p$, there exists $\be \in E$
with $\e(\be) = p'$ such that either $\a$ covers $\be$ or $\be$ covers $\a$
in $E$.}
\enddefinition

\remark{Remark \secaa.5}
\item{(i)}{It is easily checked that these are sound definitions, because 
they are invariant under isomorphism of labelled heaps.}
\item{(ii)}{Condition (a) provides a way to name the elements of $E$, which
we shall need in the sequel.
Choose a vertex $p$ of $\Gamma$.  Since $E$ is locally finite, $\e^{-1}(p)$
is a chain of $E$ isomorphic as a partially ordered set to the integers, so
one can label each element of this chain as $E(p, z)$ for some $z \in \zed$.
Adopting the convention that $E(p, x) < E(p, y)$ if $x < y$, this labelling
is unique once a distinguished vertex $E(p, 0) \in \e^{-1}(p)$ has been
chosen for each $p$.  When $E$ is understood, we will use the shorthand
$p(y)$ for $E(p, y)$.}

\head \S\secab.  Full heaps \endhead

We are now ready to recall the definition of our main object of study from
\cite{{\bf 7}}, which builds on work of Stembridge \cite{{\bf 17}} and Wildberger
\cite{{\bf 19}}.

\definition{Definition \secab.1}
Let $E$ be a fibred heap over a Dynkin diagram $\Gamma$ with generalized
Cartan matrix $A$.  If every open interval $(\a, \be)$ of $E$
such that $\e(\a) = \e(\be) = p$ and $(\a, \be) \cap \e^{-1}(p) = \emptyset$
satisfies $\sum_{\g \in (\a, \be)} a_{p, \e(\g)} = -2$, we call $E$ a
{\it full} heap.
\enddefinition

The definition says that 
either (a) $(\a, \be)$ contains precisely two elements with labels adjacent 
(via simple edges) 
to $p$, or that (b) $(\a, \be)$ contains precisely one element with label 
($q$, say) adjacent to $p$ such that there is a double edge with 
an arrow from $q$ to $p$ in the Dynkin diagram.

\definition{Definition \secab.2}
Let $R^+$ be the set of all functions $P \ra \zed^{\geq 0}$.
If $F$ is a finite labelled heap over $\Gamma$, then 
we define the {\it character}, $\cha(F)$ of $F$ to
be the element of $R^+$ such that $\cha(F)(p)$
is the number of elements of $F$ with $\e$-value $p$.  
If $\a \in R^+$, we write $\L_\a(E)$ to be the set of all convex
subheaps $F$ of $E$ with $\cha(F) = \a$.  If $F$ consists
of a single element $\a$ with $\e(\a) = p$, we will write $\cha(F) = p$ for
short, so that $\L_p(E)$ is identified with the elements of $E$ labelled
by $p$.
Since the function $\cha$ is an invariant of labelled heaps, we can extend
the definition to apply to finite heaps of $\Gamma$.
\enddefinition

\definition{Definition \secab.3}
Let $E$ be a full heap over a graph $\Gamma$ and let $k$ be a field.
Let $\Be$ be the set of proper ideals of $E$, so that $\Be$ has the structure
of a distributive lattice 
with meet and join operations $I \wedge J = I \cap J$ and
$I \vee J = I \cup J$; these operations are defined by
\cite{{\bf 7}, Lemma 2.1 (ii)}.
Let $V_E$ be the $k$-span of the set $\{v_I : I \in \Be\}$.
For any proper ideal and any finite convex subheap $L \leq E$, we
write $L \succ I$ to mean that both $I \cup L$ is an ideal and $I \cap L = 
\emptyset$, and we write $L \prec I$ to mean
that both $L \leq I$ and $I \backslash L$ is an ideal.
We define linear operators 
$X_L$, $Y_L$ and
$H_L$ on $V_E$ as follows: $$\eqalign{
X_L(v_I) &= \cases
v_{I \cup L} & \text{ if } L \succ I,\cr
0 & \text{ otherwise,}\cr
\endcases\cr
Y_L(v_I) &= \cases
v_{I \backslash L} & \text{ if } L \prec I,\cr
0 & \text{ otherwise,}\cr
\endcases\cr
H_L(v_I) &= \cases
v_I & \text{ if } L \prec I \text{ and } L \not\succ I,\cr
-v_I & \text{ if } L \succ I \text{ and } L \not\prec I,\cr
0 & \text{ otherwise.}\cr
\endcases\cr
}$$  
If $p$ is a vertex of $\Gamma$, we write $X_p$ for the linear operator
on $V_E$ given by $$\sum_{L \in \L_p(E)} X_L,$$ and we define $Y_p$ and $H_p$
similarly.  Note that although these sums are infinite, it follows from
the definitions of fibred and full heaps that at most one of the
terms in each case may act in a nonzero way on any given $v_I$.  In this
situation, we also write $p \succ I$ to mean that $L \succ I$ for some
(necessarily unique) $L \in \L_p(E)$, and analogously we write
$p \prec I$ with the obvious meaning.  Note that it is not possible for
both $p \prec I$ and $p \succ I$, because $I$ cannot contain a convex
chain $\a < \be$ with $\e(\a) = \e(\be) = p$.
\enddefinition

By \cite{{\bf 7}, Lemma 2.1}, the operators $X_L$, $Y_L$ and $H_L$ are nonzero
and well defined (see also \cite{{\bf 7}, Definition 2.6}).

\example{Example \secab.4}
The operators $X_p$, $Y_p$ and $H_p$ are of key importance in the application
of full heaps to affine Kac--Moody algebras: they respectively represent 
the action of the Chevalley generators $e_p$, $f_p$ and $h_p$ of the derived
algebra $\fg'(A)$.

For an example of these operators, consider, 
the heap $E$ shown in Figure \secea.2 in \S\secea, and let
$F$ be the ideal $\lan 2(1), 5(1) \ran$ of $E$.  Then $F$ consists of
the elements $$
\{0(y), 1(y), 4(y) : y \leq 0\} 
\cup 
\{2(y), 3(y), 5(y) : y \leq 1\} 
.$$  In this case, we have $X_4(v_F) = v_{F \cup \{4(1)\}}$, 
$X_0(v_F) = v_{F \cup \{0(1)\}}$, and $X_p(v_F) = 0$ for 
$p \in \{1, 2, 3, 5\}$.  We also have $Y_2(v_F) = v_{F \backslash \{2(1)\}}$,
$Y_5(v_F) = v_{F \backslash \{5(1)\}}$, and $Y_p(v_F) = 0$ for $p \in \{0, 1, 
3, 4\}$.  We have $$
H_p(v_F) = \cases
v_F & \text{ if } p \in \{2, 5\},\cr
-v_F & \text{ if } p \in \{0, 4\},\cr
0 & \text{ if } p \in \{1, 3\}.\cr
\endcases$$\endexample

\definition{Definition \secab.5}
Let $A$ be a simply laced generalized Cartan matrix of affine type 
and let $\Gamma$ be the corresponding Dynkin diagram, and suppose that $\mu$
is a nonidentity graph automorphism of $\Gamma$.  
Assume furthermore that (a) $\mu$ has order precisely $2$ and (b) for any 
vertex $p$, $\mu(p)$ and $p$ are not distinct adjacent vertices.
The group $\{1, \mu\}$ acts on the Dynkin
diagram $\Gamma$, and we denote the orbit containing the vertex $p$ by
$f(p) = \bar{p}$.  

The Dynkin diagram $\overline{\Gamma}$ for $\overline{\fg}$ has vertices
labelled by the orbits $\bar{p}$, and is such that if $p$ and $q$ are
distinct vertices of $\Gamma$, then $p$ and $q$
are adjacent in $\Gamma$ if and only if the (distinct) vertices $\bar{p}$
and $\bar{q}$ are adjacent in $\overline{\Gamma}$.  If $\Gamma$ contains
three vertices $p$, $\mu(p)$ and $q$ such that $q$ is adjacent to both
$p$ and $\mu(p)$, then we join $\bar{p}$ and $\bar{q}$ in $\overline{\Gamma}$
by a double edge with
an arrow pointing towards $\bar{p}$.  (It is possible for this procedure to
result in a double edge with two arrows in opposite directions.)

In the above situation, we will say that $A$ (respectively, 
$\Gamma$) {\it folds} to $\overline{A}$ (respectively, $\overline{\Gamma}$) 
via $\mu$.
\enddefinition

\definition{Definition \secab.6}
Let $A$ be a generalized Cartan matrix of affine type
with Dynkin diagram $\Gamma$, and let $E = (E, \leq, \e)$ be a full heap
over $\Gamma$.
\item{(i)}{The heap $E$ is by definition a {\it simply folded full heap}
over $\Gamma$.}
\item{(ii)}
{Suppose  there exists a diagram automorphism $\mu$ as in Definition
\secab.5 such that $A$ and $\Gamma$ fold to $\overline{A}$ and 
$\overline{\Gamma}$ respectively via $\mu$.
Suppose furthermore that whenever we have vertices $p, q$
of $\Gamma$ satisfying (a) $\mu(p) \ \C \  q$, (b) $\a \in \e^{-1}(p)$ and 
(c) $\be \in \e^{-1}(q)$, then $\a$ and $\be$ are comparable in $E$.
Then we say that
$\overline{E} := (E, \leq, f \circ \e)$ is a {\it simply folded full heap} 
over $\overline{\Gamma}$.}
\enddefinition

It is not immediate that the heap $\overline{E}$ is well defined, but this
is follows from \cite{{\bf 7}, Proposition 6.1}.  

\remark{Remark \secab.7}
All examples of full heaps in this paper will be simply folded.
\endremark

\head \S\secba. The Weyl group and skew ideals \endhead

We define the Weyl group, $W(\Gamma)$, associated to $\Gamma$ to be the
group with generators $S(\Gamma) = \{s_i \in I\}$ indexed by the vertices 
of $\Gamma$ and defining relations $$\eqalign{
s_i^2 &= 1 \text{ for all } i \in I,\cr
s_i s_j &= s_j s_i \text{ if } a_{ij} = 0,\cr
s_i s_j s_i &= s_j s_i s_j \text{ if } a_{ij} < 0 \text{ and } a_{ij} 
a_{ji} = 1,\cr
s_i s_j s_i s_j &= s_j s_i s_j s_i \text{ if } a_{ij} < 0 \text{ and } 
a_{ij} a_{ji} = 2.\cr
}$$  Note that no relation is added in the case where $a_{ij} < 0$ and $a_{ij}
a_{ji} = 4$.

\example{Example \secba.1}
Define two generalized Cartan matrices $$
A_1 = \left( \matrix
2 & -1 \cr
-2 & 2 \cr
\endmatrix \right) 
\text{\ and\ }
A_2 = \left( \matrix
2 & -2 \cr
-2 & 2 \cr
\endmatrix \right) 
.$$  Then the Weyl group corresponding to $A_1$ is $$
\lan s_1, s_2 : s_1^2 = s_2^2 = 1, (s_1 s_2)^4 = 1 \ran
,$$ isomorphic to the dihedral group of order $8$, and the Weyl group
corresponding to $A_2$ is the infinite group $$
\lan s_1, s_2 : s_1^2 = s_2^2 = 1 \ran
.$$
\endexample

\definition{Definition \secba.2 \cite{{\bf 7}, Definition 8.6}}
Let $A$ be a generalized Cartan matrix 
with Dynkin diagram $\Gamma$, and let $E$ be a simply folded full heap over
$\Gamma$.
For each vertex $i$ of $\Gamma$, we define a linear operator $S_i$ on $V_E$
by requiring that $$
S_i(v_I) = \cases
Y_i(v_I) & \quad \text{if } Y_i(v_I) \ne 0,\cr
X_i(v_I) & \quad \text{if } X_i(v_I) \ne 0,\cr
v_I & \quad \text{otherwise.}
\endcases$$  It follows that $S_i(v_I) = v_{I'}$ for some proper ideal $I'$
of $E$, so we also write $S_i(I) = I'$.
\enddefinition

The next result, whose proof is immediate from the definitions, will be useful
in the sequel.

\proclaim{Lemma \secba.3}
Maintain the notation of Definition \secba.2, and define the integer $c$
by $H_i(v_I) = c v_I$.  Then we have $$
S_i(v_I) = \cases
Y_i(v_I) & \quad \text{if } c = 1,\cr
X_i(v_I) & \quad \text{if } c = -1,\cr
v_I & \quad \text{if } c = 0.\qed\cr
\endcases$$\endproclaim

The following results show how the Weyl group acts naturally on the 
lattice $\Be$.
The main purpose of this paper is to understand this action, which remarkably
turns out to be faithful in all the well-understood cases.

\proclaim{Proposition \secba.4}
Let $A$ be a generalized Cartan matrix 
with Dynkin diagram $\Gamma$, and let $E$ be a simply folded full heap over
$\Gamma$, and let $E^*$ be the dual heap.  Let $\Be$ and $\Be^*$ be the
lattices of ideals of $E$ and $E^*$, respectively.

\item{\rm (i)}
{The Weyl group acts transitively on $\Be$ via $
s_i . b = S_i(b)
$ for each generator $s_i$ of $W$.}
\item{\rm (ii)}
{The map $* : E \ra E^*$ induces a map $*_{\Be} : \Be \ra \Be^*$ given by $$
*_{\Be}(I) = *(E \backslash I)
,$$ and this is an isomorphism of $W$-sets.}
\endproclaim

\demo{Proof}
In \cite{{\bf 7}, Proposition 8.7}, it is shown that the assignment
$s_i \mapsto -S_i$ defines a unique cyclic $kW$-module structure on $V_E$.  
The well-definedness of the action follows from this by twisting by sign, 
together with the fact that $S_i$ induces a function from $\Be$ to $\Be$.
Transitivity now follows from the cyclic module structure, completing the
proof of (i).

Using the definition of proper ideal, it is routine to check 
that the map $*_{\Be}$ of (ii) is a defined and bijective.
Given $I \in \Be$ and $I^* = *_{\Be}(I) \in \Be^*$ and a vertex $i$ of
$\Gamma$, we find that $L \succ I$ if and only if $L^* \prec^* I^*$,
and $L \prec I$ if and only if $L^* \succ^* I^*$ (where $L^* = *_{\Be}(L)$,
and $\prec^*$ and $\succ^*$ have the obvious meanings).  
From the symmetry between $Y_i$ and $X_i$ in the definition
of $S_i$ (Definition \secba.2), we can now deduce that the map
$*_{\Be}$ intertwines the two $W$-actions, completing the proof of (ii).
\qed\enddemo

The key to understanding the Weyl group action will turn out to be the root
system associated to $W$ and the generalized Cartan matrix, which we now
introduce.

Let $\Pi = \{\a_i : i\in I\}$ and let $\Pi^\vee = 
\{\a_i^\vee : i \in I\}$.  We have a $\zed$-bilinear pairing $\zed\Pi \times
\zed\Pi^\vee \ra \zed$ defined by $$
\lrt{\a_j}{\a_i^\vee} = a_{ij}
,$$ where $(a_{ij})$ is the generalized Cartan matrix.
If $k$ is a field, we extend this to a $k$-bilinear pairing by extension
of scalars.
If $v = \sum_{i \in I} \l_i \a_i$, we write $v \geq 0$ to mean that 
$\l_i \geq 0$
for all $i$, and we write $v > 0$ to mean that $\l_i > 0$ for all $i$.
We view $V = k\Pi$ as the
underlying space of a reflection representation of $W$, determined
by the equalities $s_i(v) = v - \lrt{v}{\a_i^\vee} \a_i$ for all
$i \in I$.

Indecomposable generalized 
Cartan matrices come in three mutually exclusive types
(defined in \cite{{\bf 10}, Theorem 4.3}) called {\it finite}, {\it affine} and 
{\it indefinite}.  This paper is mostly concerned with finite and
affine generalized Cartan matrices; the classification of these matrices
may be found in \cite{{\bf 10}, \S4.8}.

Following \cite{{\bf 10}, \S5}, we define a {\it real root} to be a vector
of the form $w(\a_i)$, where $w \in W$ and $\a_i$ is a basis vector.  If
$A$ is of finite type, all roots are real.  If $A$ is of affine type, there
is a unique vector $\d = \sum a_i \a_i$ such that $A\d = 0$ and the $a_i$
are relatively prime positive integers.  Although the notion of {\it imaginary
root} can be defined in general, in the affine type case the imaginary
roots are easily characterized as precisely those vectors of the form 
$n\d$ where $n$ is a nonzero integer.

A {\it root} is by definition a real or imaginary root.  
We denote the set of roots by $\Delta$, as in \cite{{\bf 10}}.
We say a root $\a$ is
positive (respectively, negative) if $\a > 0$ (respectively, $\a < 0$).
If $\a$ is a root, then so is $-\a$, and every root is either positive or
negative.  
Following \cite{{\bf 10}}, we use the symbols $\Delta, \rreal$
and $\Delta^{im}$ for the roots, real roots, and imaginary roots respectively.
We denote the positive and negative roots by $\Delta_+$ and $\Delta_-$
respectively, and we use the notation $\rreal_\pm$ and $\Delta^{im}_\pm$
with the obvious meanings.
We will also identify the positive (real and
imaginary) roots with elements of $R^+$ as in Definition \secab.2 so that
$\sum a_i \a_i$ corresponds to the function sending each $i$ to $a_i$.

\proclaim{Lemma \secba.5}
Let $E$ be a simply folded full heap over a finite graph $\Gamma$, and let
$F$ and $F'$ be proper ideals of $E$.  Then $F \backslash F'$ and
$F' \backslash F$ are finite heaps.
\endproclaim

\demo{Proof}
By \cite{{\bf 7}, Lemma 2.1 (ii)}, $F \cap F'$ is a proper ideal of $E$.  Since
$\Gamma$ is finite, \cite{{\bf 7}, Lemma 2.1 (vi)} applied to the fact that
$F \cap F' \subseteq F$ shows that $F \backslash (F \cap F')$ is finite,
so that $F \backslash F'$ is finite, as required.  The other assertion follows
similarly.
\qed\enddemo

Lemma \secba.5 ensures that the following definition makes sense.

\definition{Definition \secba.6}
Let $E$ be a simply folded full heap over a finite graph $\Gamma$, and let
$F$ and $F'$ be proper ideals of $E$.  We define the character,
$\cha(F, F') = \cha((F, F'))$,
of the ordered pair $(F, F')$ to be the function $P \ra \zed$ given by $$
\cha(F, F')(p) = \cha(F' \backslash F)(p) - \cha(F \backslash F')(p)
.$$  If $\cha(F, F') \in \rreal$,
then the proper ideals $F$ and $F'$ are said to be {\it skew}.  We define
$\Sigma(E)$ to be the set of pairs $(F, F')$ of skew proper ideals of $E$.
\enddefinition

\remark{Remark \secba.7}
Note that the relation ``is skew to'' is irreflexive and symmetric.  
The reason for the term ``skew'' will become clear when we study type $E_6$
in \S\sectd.
\endremark

\proclaim{Lemma \secba.8}
Let $E$ be a simply folded full heap over a finite graph $\Gamma$, and let
$F$ and $F'$ be proper ideals of $E$.
\item{\rm (i)}
{We have $F \subseteq F'$ if and only if $\cha(F, F')(p) \geq 0$ for all
vertices $p$ of $\Gamma$.}
\item{\rm (ii)}
{We have $F' \subseteq F$ if and only if $\cha(F, F')(p) \leq 0$ for all
vertices $p$ of $\Gamma$.}
\item{\rm (iii)}
{If $F, F', F''$ are proper ideals of $E$, then $\cha(F, F'') = \cha(F, F')
+ \cha(F', F'')$.}
\endproclaim

\demo{Proof}
To prove (i), first note that if $F \subseteq F'$, then we have 
$\cha(F \backslash F')(p) = 0$ for all vertices $p$ of $\Gamma$, and 
the assertion follows from the definition of $\cha(F, F')$.  

Conversely, suppose that $\cha(F, F')(p) \geq 0$ for all vertices
$p$ of $\Gamma$, and let $F_p = F \cap \e^{-1}(p)$ and $F'_p = 
F' \cap \e^{-1}(p)$.
We observe that $$
\cha(F, F')(p) = \cha(F'_p \backslash F_p)(p) - \cha(F_p \backslash F'_p)(p)
.$$  Since $F$ and $F'$ are proper ideals, it must be the case that one of the
sets $F_p$ and $F'_p$ is included 
in the other, which means that at least one of $\cha(F'_p \backslash F_p)(p)$
and $\cha(F_p \backslash F'_p)(p)$ is zero.  The hypothesis 
$\cha(F, F')(p) \geq 0$ forces $\cha(F_p \backslash F'_p)(p) = 0$ and
$F_p \subseteq F'_p$.  Since this is true for all $p$, we have $F \subseteq
F'$, completing the proof of (i).
The proof of (ii) follows by a similar argument.

Part (iii) also follows by using the techniques of the above paragraph.
The chains $F_p, F'_p, F''_p$ are totally ordered by inclusion and we
obtain $$
\cha(F, F'')(p) = \cha(F, F')(p) + \cha(F', F'')(p)
,$$ from which the assertion follows.
\qed\enddemo

\proclaim{Lemma \secba.9}
Let $E$ be a full heap over a graph $\Gamma$, let $\Be$ be the associated
distributive lattice and let $p$ and $q$ be vertices of
$\Gamma$ (allowing the possibility $p = q$).  
\item{\rm (i)}
{We have $
H_p X_q - X_q H_p = a_{pq} X_q
.$}
\item{\rm (ii)}{Let $F, F' \in \Be$.  Then there exists a finite sequence
$F = F_0, F_1, \ldots, F_r = F'$ of elements of $\Be$ such that for each
$0 \leq i < r$, we have either $F_{i+1} = F_i \cup \{\be_i\}$ or $F_i = 
F_{i+1} \cup \{\be_i\}$, for some element $\be_i \in E$.  Furthermore, we 
have $$
F' = S_{m_1} S_{m_2} \cdots S_{m_r} (F)
,$$ where each $m_i$ is such that $\a_{m_i}$ appears with nonzero
coefficient in $\cha(F, F')$.}
\item{\rm (iii)}{If $F, F' \in \Be$ are such that $\cha(F, F') = \pm\a_q$, 
then $\cha(s_p(F), s_p(F')) = s_p(\pm\a_q)$.}
\endproclaim

\demo{Proof}
Part (i) is \cite{{\bf 7}, Lemma 2.7 (2)}.

For part (ii), we let $F'' = F \wedge F'$.  By Lemma \secba.5, 
$F \backslash F''$
is finite, so there is a finite sequence of operators $Y_i$ such that $$
v_{F''} = Y_{i_1} Y_{i_2} \cdots Y_{i_k} v_F
,$$ where each $\a_{i_l}$ appears with nonzero coefficient in 
$\cha(F'', F)$.  Similarly, there is a sequence of
$X_i$ such that $$
v_{F'} = X_{j_1} X_{j_2} \cdots X_{j_{k'}} v_{F''}
,$$ where each $\a_{j_l}$ appears with nonzero coefficient in $\cha(F', F'')$.
It follows from the definition of $S_i$
that $$
F' = S_{j_1} S_{j_2} \cdots S_{j_{k'}} S_{i_1} S_{i_2} \cdots S_{i_k} (F)
.$$  Since $F'' = F \wedge F'$, an argument like that used to prove 
Lemma \secba.8 (i) and (ii) shows that the supports of $\cha(F', F'')$ 
and $\cha(F'', F)$ are disjoint (in the notation of Lemma \secba.8, we have
$F''_p = F_p$ if $F_p \subseteq F'_p$, and $F''_p = F'_p$ otherwise).
It then follows that each $\a_{i_l}$ and
$\a_{j_l}$ mentioned above must appear with nonzero coefficient in
$\cha(F', F'') + \cha(F'', F)$, which is equal to $\cha(F', F)$ by Lemma
\secba.8 (iii).  This completes the proof of (ii).

To prove (iii), suppose that $F, F'$ are as in the statement.  We only
need to consider the case $\cha(F, F') = \a_q$, because the other case
follows by exchanging the roles of $F$ and $F'$.
We then need to prove that $$
\cha(s_p(F), s_p(F')) = \cha(F, F') -  \lrt{\a_q}{\a_p^\vee} \a_p
$$  By Lemma \secba.9 (i), we
have $$
H_p \circ X_q - X_q \circ H_p = \lrt{\a_q}{\a_p^\vee} X_q
.$$  Now $X_q . v_F = v_{F'}$, and by definition of $H_i$, there exist 
integers $c_1, c_2 \in \{-1, 0, 1\}$ such that $H_i \circ X_q . v_F = 
c_1 v_{F'}$ and $X_q \circ H_i . v_F = c_2 v_{F'}$ and $c_1 - c_2 = 
 \lrt{\a_q}{\a_p^\vee}$.

Let us first consider the case where $\lrt{\a_q}{\a_p^\vee} = 2$, which
implies that $c_1 = 1$ and $c_2 = -1$.  The fact that 
$H_p \circ X_q . v_F = 
v_{F'}$ means that $H_p . v_{F'} = v_{F'}$, and, by Lemma \secba.3, $S_p . v_F
= Y_p . v_F$, so that $S_p$ removes an element $\g_1 \in E$ 
(with $\e(\g_1) = p$) 
from $F'$.  A similar argument shows that $S_p . v_F = X_p . v_F$, so that
$S_p$ adds an element $\g_2 \in E$ (with $\e(\g_2) = p$) to $F$.  If $\g_1 \ne 
\g_2$, it follows that $s_p(F) \subset s_p(F')$ and $$
s_p(F') \backslash s_p(F) = F' \backslash (F \cup \{\g_1\} \cup \{\g_2\})
,$$ meaning that $\cha(s_p(F), s_p(F')) = \cha(F, F') - 2\a_p$, as required.
The other possibility is that $\g_1 = \g_2$ is the unique element of
$F' \backslash F$.  It
follows that $s_p$ exchanges $F'$ and $F$, so that 
$\cha(F, F') = \a_p$ and 
$\cha(s_i(F), s_i(F')) = -\a_p = s_p(\a_p)$, as required.

The cases where $\lrt{\a_q}{\a_p^\vee} \in \{-2, -1, -0, 1\}$ follow very
similar lines, but they are simpler in the sense that situations analogous
to the $\g_1 = \g_2$ case above do not occur, and we always have 
$s_p(F) \subset s_p(F')$.  The only real change needed to the above argument 
is that the value of
$\lrt{\a_q}{\a_p^\vee}$ may lead to more than one possibility for $c_1$ (and
$c_2$).
\qed\enddemo

\head \S\secbb. Main results \endhead

In \S\secbb, we develop the main theoretical results of the paper; the
remaining sections will be devoted to the study of specific examples.

\proclaim{Theorem \secbb.1}
Let $E$ be a simply folded full heap over a (finite) Dynkin diagram 
$\Gamma$ for an affine Kac--Moody algebra. 
\item{\rm (i)}
{Let $(F, F') \in \Be \times \Be$, and let $w \in W$.  Then we have $$
\cha(w(F), w(F')) = w(\cha(F, F'))
.$$}\item{\rm (ii)}
{Let $F, F' \in \Be$, and let $w \in W$.  Then $(F, F') \in \Sigma(E)$ if
and only if \hfill\newline $(w(F), w(F')) \in \Sigma(E)$.}
\item{\rm (iii)}
{The action of $W$ on $\Be$ is faithful.}
\endproclaim

\demo{Note}
It is necessary to use pairs of heaps rather than single heaps to express
(i).  Indeed,
even if it is given that $F \subset F'$ and $w(F) \subset w(F')$, it is
not generally the case that the finite subheap $w(F') \backslash w(F)$ 
is a function of $F' \backslash F$ and $w$.
\enddemo

\demo{Proof}
To prove (i), we first deal with the case where $w = s_p$.  Let $F, F'$
be as in the statement, and let $r$ be the integer defined in Lemma \secba.9
(ii).  The proof is by induction on $r$; the case $r = 0$ is trivial, and
the case $r = 1$ is Lemma \secba.9 (iii).  The inductive step is given
by Lemma \secba.7 (iii), thus completing the proof in the case $w = s_p$.  The
proof of (i) for general $w$ then follows by a straightforward induction.

The ``only if'' direction of part (ii) follows from the fact
\cite{{\bf 10}, \S5.1} that $w$ permutes $\rreal$, 
and the ``if'' direction holds because $w$ is invertible.
To prove (iii), let $w \in W$ be such that $w \ne 1$; we will be done
if we can show the existence of a proper ideal $F''$ such that 
$w . F'' \ne F''$.  Since $w \ne 1$, it follows from standard properties
\cite{{\bf 10}, Lemma 3.11 (b)} that
$w(\a_i) = \a < 0$ for some simple root $\a_i \in \preal$.  Let $F, F'$ be
proper ideals such that $\cha(F, F') = \a_i$; such ideals exist by
\cite{{\bf 7}, Lemma 2.1 (vii)} (which is not hard to check directly in 
this case).
By (i), we have $\cha(w(F), w(F')) = \a \ne \cha(F, F')$, so we cannot have
both $w(F) = F$ and $w(F') = F'$, establishing the existence of $F''$ as
above.
\qed\enddemo

\remark{Remark \secbb.2}
It is possible to extend the above theorem to deal with imaginary roots, but
this is not interesting from the point of view of the Weyl group, because
in the situations covered by the theorem,
the action of $W$ on the imaginary roots is trivial.
\endremark

The Weyl groups associated to affine Kac--Moody algebras, which
are the main examples of Weyl groups of interest to us, are equipped with
a distinguished generator, $s_0$.

\definition{Definition \secbb.3}
Let $A$ be a generalized Cartan matrix of affine type
with Dynkin diagram $\Gamma$ and distinguished vertex $0$, and let $E$ be 
a full heap over $\Gamma$.  Let $W$ be the Weyl group associated to
$\Gamma$, and let $W_0$ be the (finite) subgroup generated by 
$S(\Gamma) \backslash \{ s_0 \}$.
If $F$ is a proper ideal of $E$, we define the {\it height}, $h(F)$,
of $F$ to be the maximal integer $t$ such that $E(0, t) \in F$.
\enddefinition

\proclaim{Lemma \secbb.4}
Maintain the notation of Definition \secbb.3, and fix $t \in \zed$.
\item{\rm (i)}{If $F, F'$ are proper ideals of $E$, then $
h(F \wedge F') = \min(h(F), h(F'))
$ and $
h(F \vee F') = \max(h(F), h(F'))
.$}
\item{\rm (ii)}{The subset $\Be_t := \{ F \in \Be : h(F) = t \}$ 
is a sublattice of $\Be$.}
\item{\rm (iii)}{If $I$ and $I'$ are proper ideals of $E$ and 
$X_i . v_I = v_{I'}$, then $h(I') = h(I) + \d_{0i}$ (the Kronecker delta).}
\item{\rm (iv)}{If $I$ and $I'$ are proper ideals of $E$ and 
$Y_i . v_I = v_{I'}$, then $h(I') = h(I) - \d_{0i}$.}
\item{\rm (v)}
{If $E$ is periodic and $\phi$ is a labelled poset automorphism of $E$,
then $$h(\phi(I)) = h(I) + k,$$ where $k \in \zed$ is such
that $\phi(E(0, 0)) = E(0, k)$.  If $k = 1$ and $t \in \zed$, then 
$\Be_t$ is a fundamental domain for the action of $\lpr$ on $\Be$.}
\item{\rm (vi)}{The group $W_0$ acts transitively on $\Be_t$, and on
the $\lpr$-orbits of $\Be$.}
\item{\rm (vii)}{The map $f : F \mapsto F \backslash
\lan E(0, 0) \ran$ defines a bijection between $\Be_0$ and the set of all
ideals (\idest including $\emptyset$ and $E_0$) of the convex subheap $E_0$ 
of $E$ given by $$
E_0 = \{
x \in E : x \not\geq E(0, 1) \text{ and } x \not\leq E(0, 0)
\}
.$$}\endproclaim

\demo{Proof}
Part (i) follows by considering the intersections of 
the chain $\e^{-1}(0)$ of $E$ with $F$, $F'$, $F \wedge F'$ and $F \vee F'$,
and (ii) is immediate from (i).  Parts (iii) and (iv) follow from the
relevant definitions.

The second assertion of (v) is a consequence of the first, and the first
assertion follows from the fact that $\phi(E(0, t)) = E(0, t + k)$ for any
$t \in \zed$, by periodicity.

To prove (vi), we first note that $W_0$ acts on $\Be_t$ by (iii) and (iv),
because none of the generators $s_i$ with $i \ne 0$ can change the height of
$F \in \Be_t$.  To prove transitivity,
choose $F, F' \in \Be_t$.  We need to prove the existence
of $w \in W_0$ such that $w(F) = F'$.  Since $h(F) = h(F')$, it follows
that $\a_0$ appears with zero coefficient in $\cha(F, F')$.  Lemma \secba.9
(ii) then completes the proof by producing the required $w$.

To prove that the heap $E_0$ in (vii) is convex, let $a, b \in E_0$ with
$a < c < b$.  Since $b \not \geq E(1, 0)$, we have $c \not\geq E(1, 0)$, and
since $a \not\leq E(0, 0)$, we have $c \not\leq E(0, 0)$.  This shows that
$c \in E_0$ and it follows that $E_0$ is convex.

Now let $I$ be the proper ideal $\lan E(0, 0) \ran$ of $E$.  Since $I$ 
is contained in every ideal of $E$ of height $0$, it follows
that the map $f : F \mapsto F \backslash I$ is a function from $\Be_0$
to the ideals of the heap $E_0$ as in the statement.  The inverse of $f$
is the map $g : g(G) = G \cup I$.  To complete the proof of (vii), it remains
to show that if $G$ is an ideal of $E_0$, then $G \cup I$ is an ideal of $E$
of height zero.
Suppose that $x \in G \cup I$ and that $y \in E$ is such that 
$y < x$.  
Since every $x \in G \cup I$ satisfies $x \not\geq E(1, 0)$ (whether
$x \in G$ or $x \in I$), we have $y \not\geq E(0, 1)$, and $G \cup I$ 
will have height $0$ if it is an ideal.
If $y \not\leq E(0, 0)$ then $y \in E_0$ and thus $y \in G$,
because $G$ is an ideal of $E_0$.  On the other hand, if $y \leq E(0, 0)$,
then $y \in I$ by definition.  In either case, $y \in G \cup I$, and we
conclude that $G \cup I$ is an ideal of $E$ (of height $0$).
\qed\enddemo

\proclaim{Theorem \secbb.5}
Maintain the notation of Definition \secbb.3, and suppose that
the heap $E$ is periodic with fundamental automorphism
$\phi$ and period $\delta$, where $\delta$ is the lowest positive imaginary
root (see Remark \secbb.6 below).
\item{\rm (i)}{For any proper ideal $F$ of $\Be$ and any $k \in \zed$, 
we have $\cha(F, \phi^k(F)) = k \d$.}
\item{\rm (ii)}{The group generated by $\phi$ is isomorphic to $\zed$ and 
acts naturally on $\Be$; if $F$ is a proper ideal of $E$, then we denote 
by $[F]$ the $\lpr$-orbit containing $F$.}
\item{\rm (iii)}{If $(F, F') \in \Sigma(E)$, then every
ideal in $[F]$ is skew to every ideal in $[F']$; in this case we will say
that the orbits $[F]$ and $[F']$ are skew.}
\item{\rm (iv)}{There is a bijection $
\w : \Be \ra \zed \times \Be_0
$ given by $\w(I) = (h(I), \phi^{-h(I)}(I))$.  The induced action of
$W$ on $\zed \times \Be_0$ is $$
w . (t, I) = (t + k, \phi^{-k}(w . I))
,$$ where $h$ is the height function on ideals and $k = h(w . I)$.  
Applying the natural projection $\zed \times \Be_0 \ra \Be_0$ and 
identifying $\Be_0$ with the $\lpr$-orbits
of $\Be$ as in Lemma \secbb.4 (v), we recover the action of $W$ on the
$\lpr$-orbits.}
\endproclaim

\demo{Proof}
Part (i) follows directly from the hypotheses.

The definition of fibred heap ensures that $\phi$ has infinite order.
If $I \in \Be$, we define $\phi(I)$ to be $\{\phi(i) : i \in I\}$; since
this is an invertible map sending proper ideals to proper ideals, (ii) follows.

Let $F$ and $F'$ be as in (iii), and let $\phi^a(F)$ and $\phi^b(F')$ be
typical elements of $[F]$ and $[F']$ respectively.  By part (i), 
Lemma \secba.8 (iii), 
and the fact that $\cha(G, G') = -\cha(G', G)$, we have $$
\cha(\phi^a(F), \phi^b(F')) = \cha(\phi^a(F), F) + \cha(F, F') + 
\cha(F', \phi^b(F')) = \cha(F, F') + (b-a)\d
.$$  By \cite{{\bf 10}, Proposition 6.3 (d)}, the sum of a real root and an 
integer multiple of $\d$ is again a real root, so $\phi^a(F)$ and 
$\phi^b(F')$ are skew, as required.

The inverse of the map given in (iv) is $\w^{-1}(t, I) = \phi^t(I)$.  The
formula for the induced action follows by a direct check.  The last assertion
is a consequence of Lemma \secbb.4 (v) and the fact that $\phi^{-k}(w.I)$
and $w.I$ are in the same orbit.
\qed\enddemo

\remark{Remark \secbb.6}
The hypotheses about the period $\d$ of $E$ used in Theorem \secbb.5 are true
for all the full heaps appearing in this paper (although not for the examples
mentioned in \S\sectf).  If the underlying Dynkin diagram of $E$ corresponds
to an untwisted affine Kac--Moody algebra, this was proved in
\cite{{\bf 7}, Lemma 7.4 (ii)}.  All the examples of full heaps in this paper
are of this type, with two exceptions, namely
examples \sectc.4 and \seceb.7, which correspond to twisted affine Kac--Moody
algebras.
For these two cases, the data in \cite{{\bf 10}, \S4.8} can be used to verify 
the hypothesis directly.
\endremark

\head \S\sectc. Permutations of $\zed$ \endhead

The first examples of Weyl group representations that we will consider
realize the affine Weyl groups as permutations of the integers, which  
Bj\"orner and Brenti \cite{{\bf 3}, p293} consider to be part of the folklore
of the subject.
In type affine $A$, the permutation representation of Example
\sectc.1 first appeared in work of
Lusztig \cite{{\bf 12}}, although without a proof of faithfulness,
and the representation was further developed by Shi \cite{{\bf 14}} and
Bj\"orner and Brenti \cite{{\bf 2}}.  In type affine $C$, the permutation
representation of Example \sectc.3 first appeared in work of B\'edard 
\cite{{\bf 1}}, again
without a proof of faithfulness, and was further studied by Shi \cite{{\bf 15}}.
These two examples are the simplest considered in this paper, in the sense
that they are the only ones for which the heap is a totally ordered set.
A unified treatment (with proofs) of affine Weyl groups of types
$A$, $B$, $C$ and $D$ as permutations of $\zed$ is given in Eriksson's
thesis \cite{{\bf 6}}.  

All the examples of full heaps in \S\sectc\  come from \cite{{\bf 7}, Appendix}.
All these heaps are periodic, and the 
dashed boxes in the diagrams indicate the repeating motif.  We use the $p(y)$
notation mentioned in Remark \secaa.5 (ii) to name individual elements in
a corresponding labelled heap.  
A dashed box in the diagram depicting a periodic heap 
will indicate the repeating motif.

\topcaption{Figure \sectc.1} The Dynkin diagram of type $A_l^{(1)} (l > 1)$
\endcaption 
\centerline{
\hbox to 2.888in{
\vbox to 1.125in{\vfill
        \includegraphics{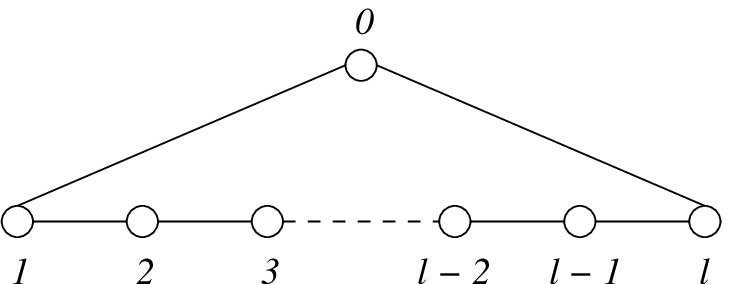}
}
\hfill}
}

\vfill\eject

\topcaption{Figure \sectc.2} A full heap, $E$, over the Dynkin diagram of type 
$A_l^{(1)} (l > 1)$
\endcaption
\centerline{
\hbox to 2.472in{
\vbox to 2.402in{\vfill
        \includegraphics{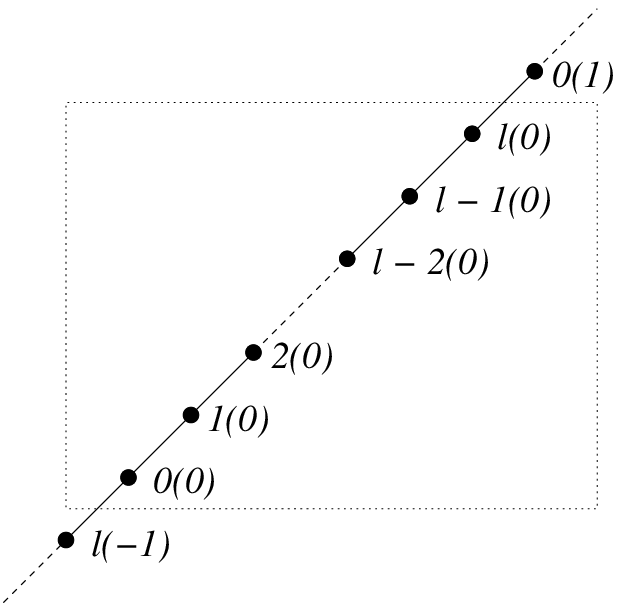}
}
\hfill}
}

\example{Example \sectc.1}
Consider the full heap shown in Figure \sectc.2 over the affine Dynkin
diagram shown in Figure \sectc.1.  In this case, every proper 
ideal is principal, so the proper ideals are precisely the set $$
\{p(y) : 0 \leq p \leq l, \ y \in \zed\}
.$$  Because the proper ideals are totally ordered, the map $\z : \Be \ra \zed$
defined by $\z(\lan p(y) \ran) = (l+1)y + p + 1$ is an isomorphism of 
totally ordered
sets (where $\zed$ is ordered in the usual way).  In this way, the action of
$W$ on $\Be$ induces an action of $w$ on $\zed$, which is faithful by
Theorem \secbb.1 (iii).

With these identifications, the action of $s_i$ on $\zed$ is as follows: $$
s_i (z) = \cases
z + 1 & \text{ if } z \equiv i \mod (l+1),\cr
z - 1 & \text{ if } z \equiv i+1 \mod (l+1),\cr
z & \text{ otherwise.}
\endcases
$$  Thus, we recover the familiar realization of the affine Weyl group of
type $A$ as permutations of the integers, as described by Lusztig 
\cite{{\bf 12}}, together with a proof that this representation 
is faithful.
\endexample

\topcaption{Figure \sectc.3} The Dynkin diagram of type $A_1^{(1)}$
\endcaption 
\centerline{
\hbox to 0.888in{
\vbox to 0.388in{\vfill
        \includegraphics{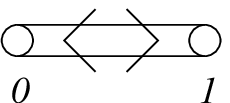}
}
\hfill}
}

\remark{Remark \sectc.2}
The case $l = 1$ of Example \sectc.1 can also be checked to give a 
faithful representation, and the analogue of the heap in Figure \sectc.2
is indeed a full heap over the Dynkin diagram shown in Figure \sectc.3.
\endremark

\topcaption{Figure \sectc.4} The Dynkin diagram of type $C_l^{(1)} (l > 1)$
\endcaption 
\centerline{
\hbox to 2.888in{
\vbox to 0.388in{\vfill
        \includegraphics{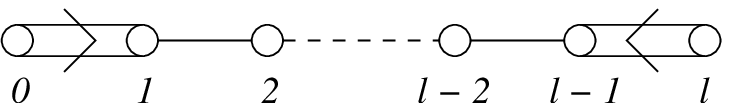}
}
\hfill}
}

\topcaption{Figure \sectc.5} A full heap, $E$, over the Dynkin 
diagram of type $C_l^{(1)} (l > 1)$
\endcaption
\centerline{
\hbox to 2.152in{
\vbox to 2.777in{\vfill
        \includegraphics{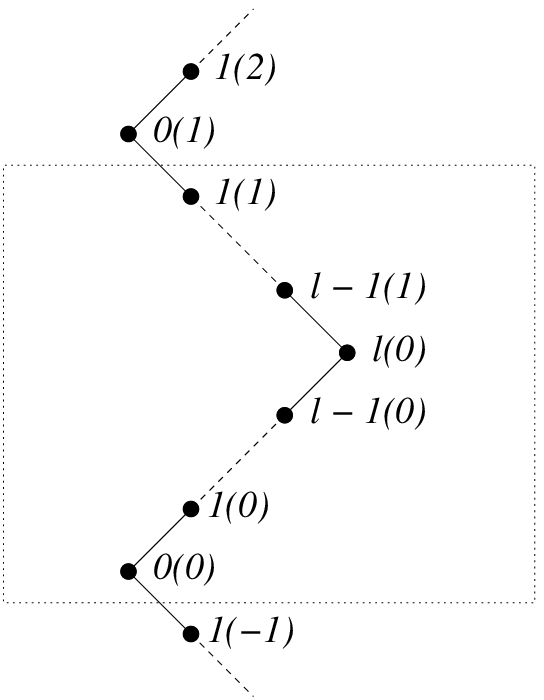}
}
\hfill}
}

\example{Example \sectc.3}
Consider the self-dual full heap shown in Figure \sectc.5 over the 
affine Dynkin
diagram shown in Figure \sectc.4.  As in Example \sectc.1, every proper ideal is 
principal, so the proper ideals are precisely the set $$
\{p(y) : 0 \leq p \leq l, \ y \in \zed\}
.$$  Because the proper ideals are totally ordered, we can again define
an isomorphism $\z : \Be \ra \zed$ of totally ordered sets, as follows: $$
\z(\lan p(y) \ran) = \cases
2ly + 1  & \text{ if } p = 0,\cr
2ly + l + 1  & \text{ if } p = l,\cr
ly + p + 1  & \text{ if } p \not\in \{0, l\} \text{ and } 
y \text{ is even},\cr
l(y+1) - p + 1  & \text{ if } p \not\in \{0, l\} \text{ and } 
y \text{ is odd}.\cr
\endcases
$$  With these identifications, the action of $s_j$ on $\zed$ is $$
s_j (z) = \cases
z + 1 & \text{ if } z \equiv \pm j \mod 2l,\cr
z - 1 & \text{ if } z \equiv (\pm j) + 1 \mod 2l,\cr
z & \text{ otherwise.}
\endcases
$$

Thus, we recover the familiar realization of the affine Weyl group of
type $C$ as permutations of the integers, as described by B\'edard 
\cite{{\bf 1}}, together with a proof that this representation 
is faithful.
\endexample

\topcaption{Figure \sectc.6} The Dynkin diagram of type $A_{2l-1}^{(2)}$
\endcaption
\centerline{
\hbox to 3.569in{
\vbox to 1.138in{\vfill
        \includegraphics{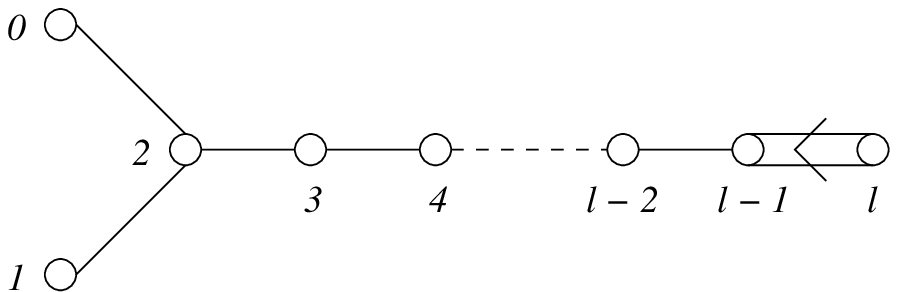}
}
\hfill}
}

\vfill\eject

\topcaption{Figure \sectc.7} A full heap, $E$, over the Dynkin 
diagram of type $A_{2l-1}^{(2)}$
\endcaption
\centerline{
\hbox to 2.152in{
\vbox to 3.527in{\vfill
        \includegraphics{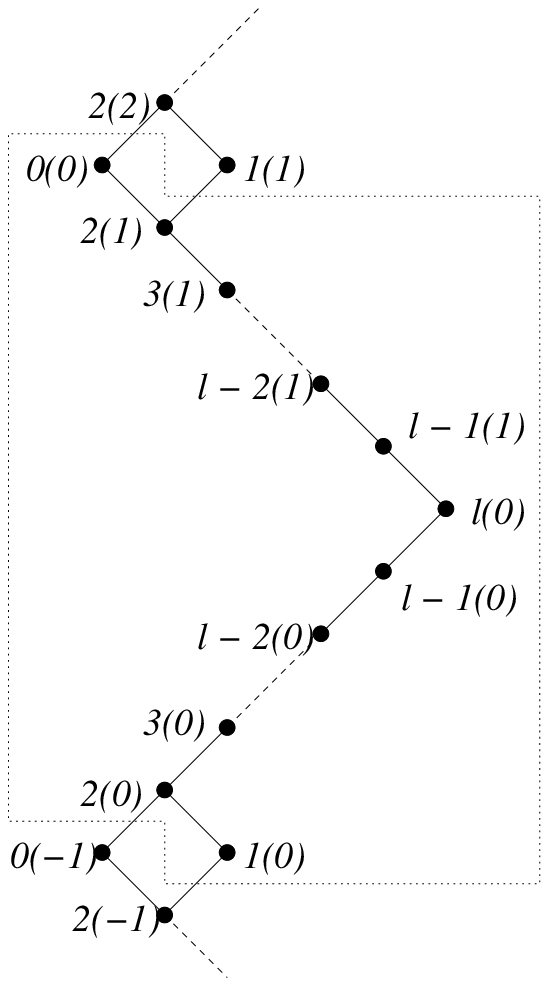}
}
\hfill}
}

\example{Example \sectc.4}
Consider the self-dual full heap $E$ of Figure \sectc.7, over the 
Dynkin diagram of type $A_{2l-1}^{(2)}$ shown in Figure \sectc.6.  
All proper ideals of $E$ are principal, 
except those of the form $\lan 0(y), 1(y+1) \ran$.  We refine the order
on the ideals to a total one by stipulating that, for all $y \in \zed$,
$\lan 0(y) \ran < \lan 1(y+1) \ran$.

With this refinement, we can define
an isomorphism $\z : \Be \ra \zed$ of totally ordered sets, as follows: $$
\eqalign{
\z(\lan 0(y-1), 1(y) \ran) &= 2ly + 2,\cr
\z(\lan p(y) \ran) &= \cases
2l(y+1) & \text{ if } p = 0,\cr
2ly + 1 & \text{ if } p = 1,\cr
ly + p + 1& \text{ if } 1 < p < l \text{ and } 
y \text{ is even},\cr
l(y+1) - p + 1  & \text{ if } 1 < p < l \text{ and } 
y \text{ is odd}.\cr
2ly + l + 1 & \text{ if } p = l.\cr
\endcases
}$$ 

With these identifications, the action of $s_j$ on $\zed$ (if 
$j \not\in \{0, 1 \}$) is $$
s_j (z) = \cases
z + 1 & \text{ if } z \equiv \pm j \mod 2l,\cr
z - 1 & \text{ if } z \equiv (\pm j) + 1 \mod 2l,\cr
z & \text{ otherwise.}
\endcases
$$  We have $$
s_0 (z) = \cases
z + 1 & \text{ if } z \equiv \pm 1 \mod 2l,\cr
z - 1 & \text{ if } z \equiv 0 \text{ or } z \equiv 2
\mod 2l,\cr 
z & \text{ otherwise,}
\endcases
$$ and $$
s_1 (z) = \cases
z + 2 & \text{ if } z \equiv -1
\text{ or } z \equiv 0 \mod 2l,\cr
z - 2 & \text{ if } z \equiv 1 \text{ or } z \equiv 2
\mod 2l,\cr 
z & \text{ otherwise.}
\endcases
$$

Thus, we recover the familiar realization of the affine Weyl group of
type $B$ as permutations of the integers, together with a proof that this
representation is faithful.
\endexample

\topcaption{Figure \sectc.8} The Dynkin diagram of type $D_l^{(1)}$
\endcaption
\centerline{
\hbox to 3.486in{
\vbox to 1.138in{\vfill
        \includegraphics{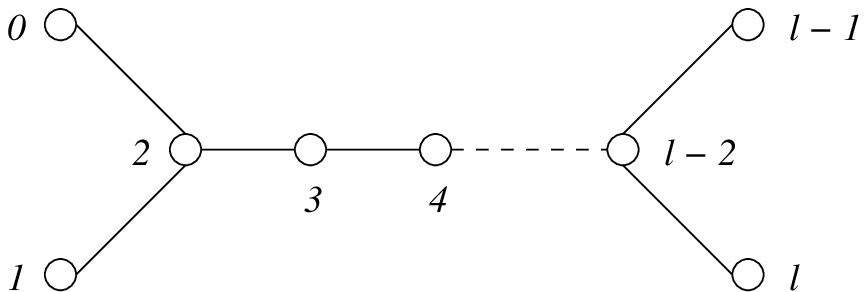}
}
\hfill}
}

\vfill\eject

\topcaption{Figure \sectc.9} A full heap, $E$, over the Dynkin diagram of type 
$D_l^{(1)}$
\endcaption
\centerline{
\hbox to 2.180in{
\vbox to 3.902in{\vfill
        \includegraphics{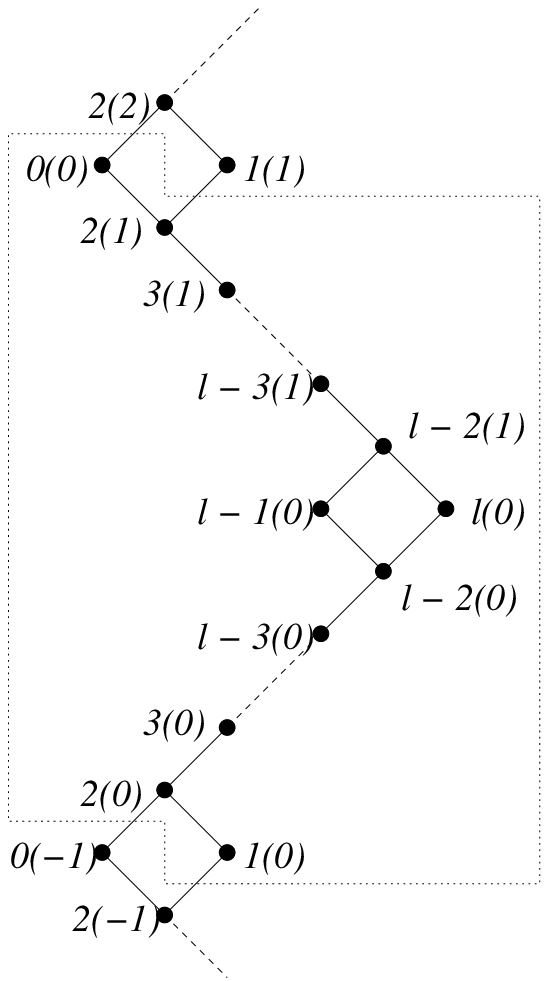}
}
\hfill}
}

\example{Example \sectc.5}
Consider the self-dual full heap shown in Figure \sectc.9 over the 
affine Dynkin 
diagram shown in Figure \sectc.8.  In this case, it is not true that every proper
ideal is principal.  The only non-principal proper ideals are 
are those of the form 
$\lan 0(y), 1(y+1) \ran$ or $\lan l-1(y), l(y) \ran$, where 
$y \in \zed$.
The set of proper ideals is not totally ordered by inclusion, but 
we may refine the order to a total one by stipulating 
that, for all $y \in \zed$,
$\lan l-1(y) \ran < \lan l(y) \ran$ and $\lan 0(y) \ran < \lan 1(y+1) \ran$.

\vfill\eject

With this refinement, we can define
an isomorphism $\z : \Be \ra \zed$ of totally ordered sets, as follows: $$
\eqalign{
\z(\lan 0(y-1), 1(y) \ran) &= 2ly + 2,\cr
\z(\lan l-1(y), l(y) \ran) &= 2ly + l + 2,\cr
\z(\lan p(y) \ran) &= \cases
2l(y+1) & \text{ if } p = 0,\cr
2ly + p & \text{ if } p = 1,\cr
2ly + p + 1 & \text{ if } p \in \{l-1, l\},\cr
ly + p + 1  & \text{ if } 1 < p < l-1 \text{ and } 
y \text{ is even},\cr
l(y+1) - p + 1  & \text{ if } 1 < p < l-1 \text{ and } 
y \text{ is odd}.\cr
\endcases
}$$ 

With these identifications, the action of $s_j$ on $\zed$ (if 
$j \not\in \{0, 1, l-1, l\}$) is $$
s_j (z) = \cases
z + 1 & \text{ if } z \equiv \pm j \mod 2l,\cr
z - 1 & \text{ if } z \equiv (\pm j) + 1 \mod 2l,\cr
z & \text{ otherwise.}
\endcases
$$  If $j \in \{0, l-1\}$, then we have $$
s_j (z) = \cases
z + 1 & \text{ if } z \equiv j - 1 + c(j) 
\text{ or } z \equiv j + 1 + c(j) \mod 2l,\cr
z - 1 & \text{ if } z \equiv j + c(j) \text{ or } z \equiv j + 2 + c(j) 
\mod 2l,\cr 
z & \text{ otherwise,}
\endcases
$$ where we define $c(0) = 0$ and $c(l-1) = 1$.  Finally, if
$j \in \{1, l\}$, we have $$
s_j (z) = \cases
z + 2 & \text{ if } z \equiv j - 2 + c(j) 
\text{ or } z \equiv j - 1 + c(j) \mod 2l,\cr
z - 2 & \text{ if } z \equiv j + c(j) \text{ or } z \equiv j + 1 + c(j) 
\mod 2l,\cr 
z & \text{ otherwise,}
\endcases
$$ where we define $c(1) = 0$ and $c(l) = 1$.

Thus, we recover the familiar realization of the affine Weyl group of
type $D$ as permutations of the integers,
together with a proof that this representation is faithful.
\endexample

By comparing the four examples above, the reader
may correctly suspect that the further the proper ideals are from being 
totally ordered by inclusion, the less helpful it is to think of the action 
of $W$ on $\Be$ as a periodic permutation of the integers, even though 
this can be done in principle.

\head \S\sectd. Geometrical examples \endhead

The application of Theorem \secbb.1 to type $E_6$ has some interesting
connections with geometry, as we shall now show by considering the full
heap $E$ in Figure \sectd.2 over the Dynkin diagram of Figure \sectd.1.  
(The heap $E$ is not self-dual, and by Proposition \secba.4 (ii),
we could equally well have started with the dual heap, $E^*$.)

\topcaption{Figure \sectd.1} The Dynkin diagram of type $E_6^{(1)}$
\endcaption
\centerline{
\hbox to 2.138in{
\vbox to 1.319in{\vfill
        \includegraphics{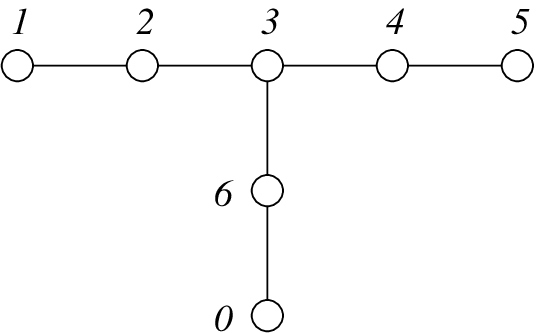}
}
\hfill}
}

\topcaption{Figure \sectd.2} A full heap, $E$, over the Dynkin diagram of 
type $E_6^{(1)}$
\endcaption
\centerline{
\hbox to 1.777in{
\vbox to 2.527in{\vfill
        \includegraphics{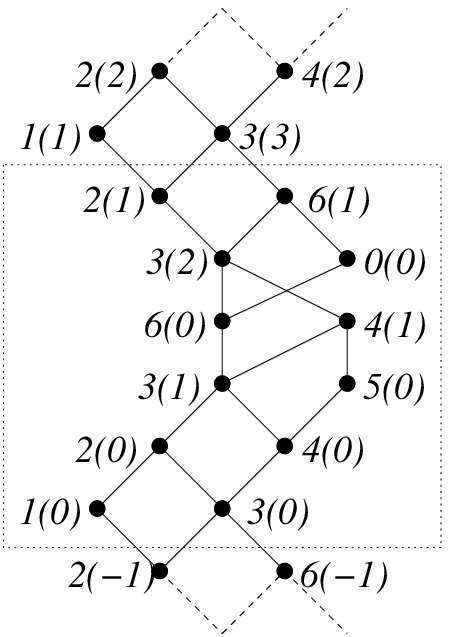}
}
\hfill}
}

The heap $E$ is periodic with fundamental automorphism $\phi$ and period
$\d = \a_1 + 2\a_2 + 3\a_3 + 2\a_4 + \a_5 + 2\a_6 + \a_0$.
By explicit enumeration, we find that there are $27$ orbits of proper ideals
of $E$ under the action of $\lpr$.

A complete set of representatives, together
with labels, is as follows: $$
\matrix
\hfill
\esla_{1} &=& \langle 1(0) \rangle,\hfill \quad & \hfill 
\esla_{2} &=& \langle 2(1) \rangle,\hfill \quad & \hfill 
\esla_{3} &=& \langle 3(2) \rangle, \hfill \cr \hfill
\esla_{4} &=& \langle 4(1), 6(0) \rangle, \hfill \quad & \hfill 
\esla_{5} &=& \langle 6(0), 5(0) \rangle, \hfill \quad & \hfill 
\esla_{6} &=& \langle 6(0) \rangle, \hfill \cr \hfill
\eslb_{1} &=& \langle 1(1), 0(0) \rangle,\hfill \quad & \hfill 
\eslb_{2} &=& \langle 2(1), 0(0) \rangle,\hfill \quad & \hfill 
\eslb_{3} &=& \langle 3(2), 0(0) \rangle, \hfill \cr \hfill
\eslb_{4} &=& \langle 4(1), 0(0) \rangle, \hfill \quad & \hfill 
\eslb_{5} &=& \langle 5(0), 0(0) \rangle, \hfill \quad & \hfill 
\eslb_{6} &=& \langle 0(0) \rangle , \hfill \cr \hfill
\eslc_{12} &=& \langle 6(1) \rangle,\hfill \quad & \hfill
\eslc_{13} &=& \langle 2(1), 6(1) \rangle,\hfill \quad & \hfill
\eslc_{14} &=& \langle 3(0) \rangle, \hfill \cr \hfill
\eslc_{15} &=& \langle 4(0) \rangle,\hfill \quad & \hfill
\eslc_{16} &=& \langle 5(0) \rangle,\hfill \quad & \hfill
\eslc_{23} &=& \langle 6(1), 1(1) \rangle, \hfill \cr \hfill
\eslc_{24} &=& \langle 1(0), 3(0) \rangle,\hfill \quad & \hfill
\eslc_{25} &=& \langle 4(0), 1(0) \rangle,\hfill \quad & \hfill
\eslc_{26} &=& \langle 5(0), 1(0) \rangle, \hfill \cr \hfill
\eslc_{34} &=& \langle 2(0) \rangle,\hfill \quad & \hfill
\eslc_{35} &=& \langle 4(0), 2(0) \rangle,\hfill \quad & \hfill
\eslc_{36} &=& \langle 5(0), 2(0) \rangle, \hfill \cr \hfill
\eslc_{45} &=& \langle 3(1) \rangle,\hfill \quad & \hfill
\eslc_{46} &=& \langle 3(1), 5(0) \rangle,\hfill \quad & \hfill
\eslc_{56} &=& \langle 4(1) \rangle. \hfill \cr
\endmatrix
$$

\proclaim{Proposition \sectd.1}
In the action of $W$ on the $\lpr$-orbits of $\Be$, the Coxeter 
generators are
represented by the following products of six transpositions: $$\eqalign{
s_1 &\mapsto (\esla_{1} \esla_{2})(\eslb_{1} \eslb_{2})(\eslc_{13} \eslc_{23})(\eslc_{14} \eslc_{24})(\eslc_{15} \eslc_{25})(\eslc_{16} \eslc_{26}), \cr
s_2 &\mapsto (\eslc_{12} \eslc_{13})(\esla_{2} \esla_{3})(\eslb_{2} \eslb_{3})(\eslc_{24} \eslc_{34})(\eslc_{25} \eslc_{35})(\eslc_{26} \eslc_{36}), \cr
s_3 &\mapsto (\eslc_{13} \eslc_{14})(\eslc_{23} \eslc_{24})(\esla_{3} \esla_{4})(\eslb_{3} \eslb_{4})(\eslc_{35} \eslc_{45})(\eslc_{36} \eslc_{46}), \cr
s_4 &\mapsto (\eslc_{14} \eslc_{15})(\eslc_{24} \eslc_{25})(\eslc_{34} \eslc_{35})(\esla_{4} \esla_{5})(\eslb_{4} \eslb_{5})(\eslc_{46} \eslc_{56}), \cr
s_5 &\mapsto (\eslc_{15} \eslc_{16})(\eslc_{25} \eslc_{26})(\eslc_{35} \eslc_{36})(\eslc_{45} \eslc_{46})(\esla_{5} \esla_{6})(\eslb_{5} \eslb_{6}), \cr
s_6 &\mapsto (\eslc_{23} \eslb_{1})(\eslc_{13} \eslb_{2})(\eslc_{12} \eslb_{3})(\esla_{4} \eslc_{56})(\esla_{5} \eslc_{46})(\esla_{6} \eslc_{45}), \cr
s_0 &\mapsto (\esla_{1} \eslb_{1})(\esla_{2} \eslb_{2})(\esla_{3} \eslb_{3})(\esla_{4} \eslb_{4})(\esla_{5} \eslb_{5})(\esla_{6} \eslb_{6}).\cr
}$$\endproclaim

\demo{Proof}
From the definition of $S_i$, it is clear that $S_i$ commutes with $\phi$, and
thus that the action of $W$ on $\Be$ commutes with the action of $\lan \phi
\ran$ on $\Be$.  The formulae for the generators may be checked by a 
routine but rather lengthy calculation.
\qed\enddemo

The significance of the above result is that the Coxeter group of type
$E_6$ (which is the subgroup $W(E_6)$ 
generated by $s_1, s_2, \ldots, s_6$) is well
known to be the automorphism group of a famous configuration of 27 lines
on a cubic surface (see \cite{{\bf 8}, Theorem V.4.9} and \cite{{\bf 8}, Exercise
V.4.11 (b)}).  We now explain how the representation
of $W(E_6)$ on the orbits of $\Be$ is isomorphic, as a permutation group, to
the action of $W(E_6)$ on the 27 lines.  Coxeter \cite{{\bf 5}, \S1} gave
the permutations representing the generators explicitly.  
In \cite{{\bf 5}}, the symbols $\esla_i, \eslb_i, \eslc_{ij}$ are the names
of the 27 lines,  and the correspondence between our notation for the 
group generators and the notation of \cite{{\bf 5}} is $$
s_1 \mapsto (1\ 2), \quad
s_2 \mapsto (2\ 3), \quad
s_3 \mapsto (3\ 4), \quad
s_4 \mapsto (4\ 5), \quad
s_5 \mapsto (5\ 6), \quad
s_6 \mapsto Q.
$$  Under these identifications, the permutations representing the action of
$W(E_6)$ in \cite{{\bf 5}, \S1} agree with those given by Proposition \sectd.1.
It follows that Proposition \sectd.1 gives an explicit action of the affine
Weyl group of type $E_6$ on the 27 lines.

We now summarise the geometric relationship between the 27 lines.
(More details may be found in \cite{{\bf 5}, \S1} or \cite{{\bf 8}, \S V.4}; note
that Hartshorne uses the notation $E_i$, $F_{ij}$, $G_i$ for $\esla_i$,
$\eslc_{ij}$, $\eslb_i$ respectively.)  Any two distinct lines that do
not intersect are skew.  The lines $\esla_1,
\ldots, \esla_6$ are mutually skew, as are the lines
$\eslb_1, \ldots, \eslb_6$.  The five skew lines $\esla_2, \ldots, \esla_6$
have a common transversal, namely $\eslb_1$, and so on.  The line
$\eslc_{ij}$ intersects $\esla_k$ (respectively, $\eslb_l$) if and only if
$k \in \{i, j\}$ (respectively, $l \in \{i, j\}$). The lines
$\eslc_{ij}$ and $\eslc_{kl}$ intersect if $\{i, j\} \cap \{k, l\} = 
\emptyset$; otherwise they are skew.  (Note that the action of $s_0$ in
Proposition \sectd.1 also preserves this relationship.)

The generators $s_1, \ldots, s_5$
act on the subscripts of the lines by transpositions, as suggested in
the correspondence of the previous paragraph.  For example, $\eslc_{24}$
is moved by $s_1$ to $\eslc_{14}$, by $s_2$ to $\eslc_{34}$, by
$s_3$ to $\eslc_{23}$ and by $s_4$ to $\eslc_{25}$.  The only orbit of
$\Be$ moved to a different $\lpr$-orbit by all
of $s_1, s_2, s_3$ and $s_4$ is the one containing $\lan 1(0), 3(0) \ran$.  
Since the action on the 27 lines is transitive (as can be checked
from the formulae in \cite{{\bf 5}, \S1}) and the action of $W_0$ on 
the $\lpr$-orbits is transitive (by Lemma \secbb.4 (vi)), we see that the
dictionary between the $\lpr$-orbits
of $\Be$ and the 27 lines is unique.

Remarkably, there is a concise description of the incidence relations 
between the 27 lines in terms of $\Be$ and the root system alone; this 
is the reason for the term ``skew'' introduced in Definition \secba.6.
The proof we give below is not conceptual, and it would be nice to know a 
reason why this should be true.

\proclaim{Proposition \sectd.2}
Let $F, F'$ be two proper ideals of the heap $E$ of Figure \sectd.2 over the Dynkin
diagram of type $E_6^{(1)}$ in Figure \sectd.1.  Let $l([F])$ and $l([F'])$ be the
corresponding lines on the cubic surface.  Then the orbits $[F]$ and $[F']$ 
are skew (in the sense of Theorem \secbb.5 (iii)) if and only if 
$l([F])$ is skew to $l([F'])$ (in the geometric sense).
\endproclaim

\demo{Proof}
The affine Weyl group $W$ preserves the heap-theoretic notion of skewness
by Theorem \secbb.1 (ii).  It also preserves the geometric notion of
skewness: the subgroup $W_0$ is known to preserve geometric skewness, and
the generator $s_0$ also preserves it, as can be checked directly from 
the incidence relations and the formula in Proposition \sectd.1.  By
transitivity of the actions, it
suffices to check the assertion for a fixed $F$.  

Let us choose $F = \lan 1(0) \ran$, so that $l(F) = \esla_1$.  In this
case, $\esla_1$ is skew to all other lines except (a) those of the form
$\eslb_i$ for $2 \leq i \leq 6$ and (b) those of the form $\eslc_{1i}$
for $2 \leq i \leq 6$.  Comparing $F$ with representatives of the orbits
corresponding to the other lines, we find that $[F]$ is skew to
all orbits except those corresponding to the aforementioned lines.
\qed\enddemo

Proposition \sectd.2 may be restated in the language of algebraic geometry
as follows, where $l_1 . l_2 \in \zed$ denotes the intersection number of 
the lines $l_1$ and $l_2$ as described in \cite{{\bf 8}, Theorem V.1.1}.

\proclaim{Corollary \sectd.3}
Let $F, F'$ be two proper ideals of the heap $E$ of Figure \sectd.2 over the Dynkin
diagram of type $E_6^{(1)}$ in Figure \sectd.1, and let $l([F])$ and $l([F'])$ be the
corresponding lines on the cubic surface.  Then the intersection number
$l([F]).l([F'])$ is given by $$
l([F]).l([F']) = \cases
-1 & \text{ if } \cha(F, F') = t\d \text{ for some } t \in \zed,\cr
0 & \text{ if } \cha(F, F') \in \rreal,\cr
1 & \text{ otherwise.}\cr
\endcases$$\endproclaim

\demo{Proof}
The conditions in the statement are easily checked (using the
definitions and \cite{{\bf 8}, Theorem V.1.1, Theorem V.4.9}) to be equivalent to
the respective conditions (a) $[F]$ and $[F']$ are equal, (b) $[F]$ 
and $[F']$ are skew and (c) $[F]$ and $[F']$ are neither skew nor equal.
\qed\enddemo

Note that an explicit description of the (faithful) action of the affine 
Weyl group $W$ on $\Be$ itself may be obtained by using 
Theorem \secbb.5 (iv).

An analogous construction to the above one for type $E_6$ can also be
performed for type $E_7$.  The Dynkin diagram of type $E_7^{(1)}$ is shown
in Figure \sectd.3, and the unique (and, therefore, self-dual) full heap 
over this graph is shown in Figure \sectd.4.

\topcaption{Figure \sectd.3} The Dynkin diagram of type $E_7^{(1)}$
\endcaption
\centerline{
\hbox to 3.138in{
\vbox to 0.819in{\vfill
        \includegraphics{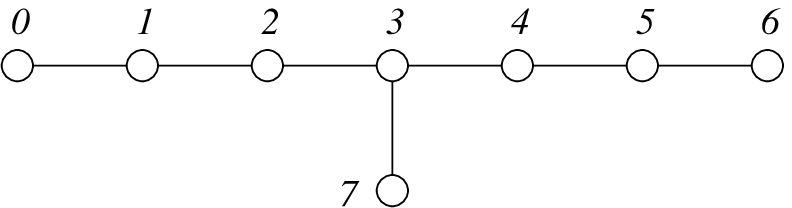}
}
\hfill}
}

\topcaption{Figure \sectd.4} A full heap over the Dynkin diagram of type 
$E_7^{(1)}$
\endcaption
\centerline{
\hbox to 2.277in{
\vbox to 2.777in{\vfill
        \includegraphics{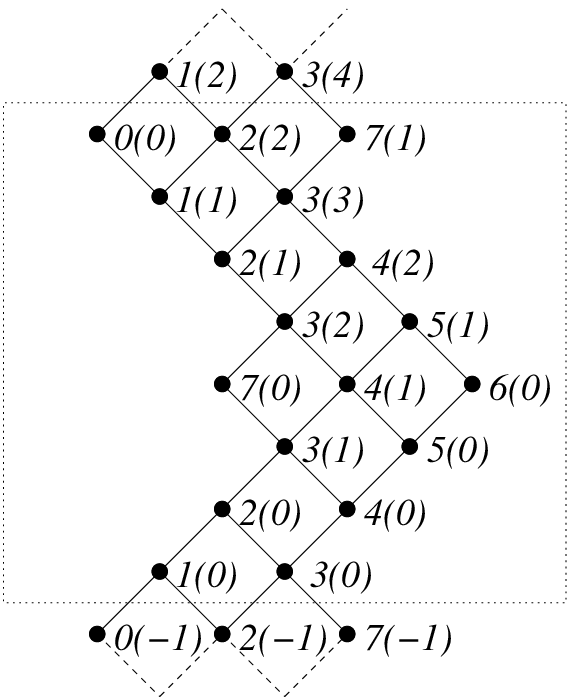}
}
\hfill}
}

\remark{Remark \sectd.4}
It is known (see for example \cite{{\bf 13}, \S4}) that the unique normal subgroup
of index 2 in the finite Weyl group of type $E_7$ is the automorphism group
of a certain configuration of lines, namely $56$ lines on the Del Pezzo 
surface of degree two defined as the double cover of the projective plane 
branched over the quartic.  Since there are $56$ orbits of proper ideals
of the heap $E$ in Figure \sectd.4 under the action of $\lpr$, we expect that
the $E_6$ construction presented here in detail, when imitated for type
$E_7$, should produce a permutation representation that is closely
related to this geometrical situation.
\endremark

\head \S\secea. The binary path representation in type affine $B$ \endhead

In this section, we will apply our theory to the full heaps corresponding
to the spin representations of the simple Lie algebra of type $B$.

Consider the self-dual full heap shown on the left of in Figure \secea.2 over 
the affine Dynkin diagram shown in Figure \secea.1.  The right hand 
side of Figure \secea.2 shows the convex subheap $E_0$ of 
Lemma \secbb.4 (vii).  

\topcaption{Figure \secea.1} The Dynkin diagram of type $B_l^{(1)}$
\endcaption
\centerline{
\hbox to 3.069in{
\vbox to 1.138in{\vfill
        \includegraphics{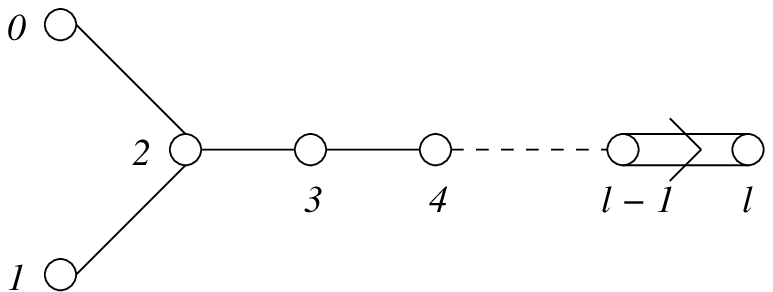}
}
\hfill}
}

\vfill\eject

\topcaption{Figure \secea.2} A full heap, $E$, over the Dynkin diagram of type 
$B_l^{(1)}$ for $l = 5$, and the convex subheap $E_0$
\endcaption
\centerline{
\hbox to 3.347in{
\vbox to 2.527in{\vfill
        \includegraphics{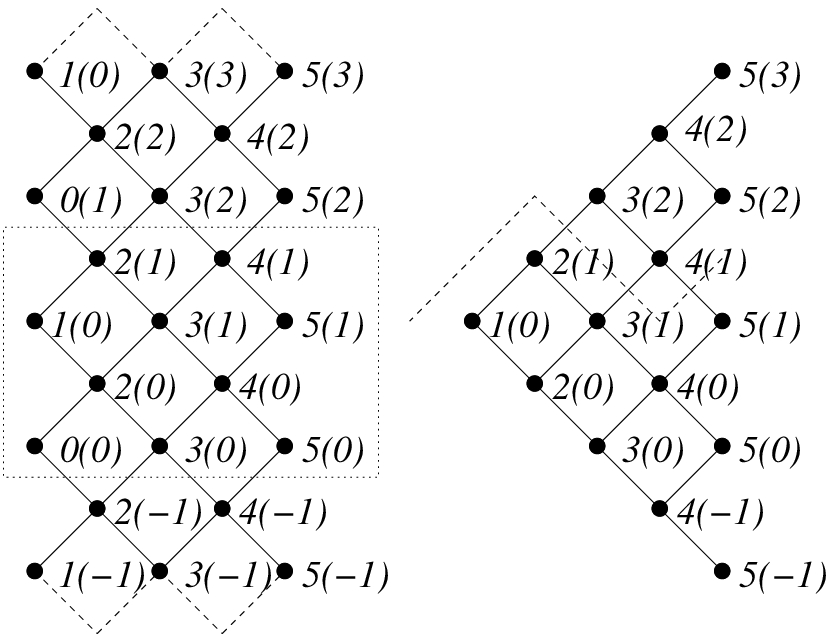}
}
\hfill}
}

In this example, 
it is convenient to regard the vertices of $E_0$ as occupying positions in
a $\zed \times \zed$ grid.  More precisely, in the heap $E_0$ arising
from type $B_l^{(1)}$, the set of occupied positions is $$
\{ (x, y) \in \zed \times \zed : 
0 \leq x < l \text{ and } -x \leq y \leq x \text{ and } x - y \in 2 \zed
\}
.$$  Note that
the smallest $y$ such that $(i, y) \in E_0$ is $y = -i$.

\definition{Definition \secea.1}
Let $E_0$ be as above.
To each ideal $G$ of $E_0$ and integer $i$ with $0 \leq i < l$, 
we define $$g(G, i) = \max \{y : (i, y) \in G\} + 1,$$ with the convention
that $g(G, i) = -i-1$ if $\{y : (i, y) \in G\}$ is empty.  
\enddefinition

\definition{Definition \secea.2}
A {\it binary path} of type $B_l^{(1)}$ is by definition a function $$
f : \{x \in \zed : -1 \leq x < l \} \ra \zed \}
$$ such that for all $0 \leq i < l$ we have $f(-1) = 0$ and 
$f(i+1) = f(i) \pm 1$.  We denote the set of all binary paths of type
$B_l^{(1)}$ by $P(B_l^{(1)})$.
\enddefinition

\proclaim{Lemma \secea.3}
There is a bijection $\pi$ from the set of binary paths of type $B_l^{(1)}$
to the set of ideals of the heap $E_0$ arising from type $B_l^{(1)}$, where
we regard an ideal of $E_0$ as a subset of $\zed \times \zed$ as above.
Explicitly, we have $$
\pi(f) = \{(x, y) \in E_0 : y < f(x)\}
.$$\endproclaim

\demo{Proof}
For $-1 \leq j < l$, define $$
(\iota(G))(j) = \cases
0 & \text{ if } j = -1,\cr
g(G, j) & \text{ otherwise,}
\endcases
$$ where the function $g$ is as in Definition \secea.1.  Note that 
if $(x, y) \in G$, then both of $(x-1, y-1)$ and $(x+1, y-1)$ must lie in $G$
because $G$ is an ideal; this is clear from the nature of the covering
relations in the Hasse diagram of $E_0$, shown in Figure \secea.2.  We also
have $g(G, 0) = \pm 1$.  Using these
observations and the definitions, we find that $\iota(G)$ is a binary path of
type $B_l^{(1)}$.  The same techniques show that $\pi(f)$ (as in the statement)
will be an ideal of $E_0$.  A routine check shows that $\pi(\iota(G)) = 
\iota(\pi(G)) = G$, and the assertions now follow.
\qed\enddemo

\example{Example \secea.4}
The dashed line in the depiction of the heap $E_0$ in Figure \secea.2 connects
the points in the binary path $f$ for which $f(-1) = f(3) = 0$, $f(0) = 
f(2) = f(4) = 1$ and $f(1) = 2$.  The corresponding ideal of $E_0$ is
$G = \lan 2(1), 5(1) \ran$, which contains $10$ elements of $E_0$.
The corresponding ideal of $E$ is $\lan 2(1), 5(1), 0(0) \ran = 
\lan 2(1), 5(1) \ran$.
\endexample

\proclaim{Corollary \secea.5}
The heap $E_0$ arising from type $B_l^{(1)}$ has $2^l$ ideals (including
$\emptyset$ and $E_0$).
\endproclaim

\demo{Proof}
This follows from Lemma \secea.3 and the observation that there 
are $2^l$ binary paths of type $B_l^{(1)}$.
\qed\enddemo

It is convenient to represent binary paths of type $B_l^{(1)}$ by the
set $\pml$ of functions
from $\{1, 2, \ldots, l\}$ to $\{+, -\}$, which we will denote by strings
of length $l$ from the alphabet $\{+, -\}$.  These also index the ideals
of $E_0$ by Lemma \secea.3.  

\definition{Definition \secea.6}
Given a binary path $f$, we define
the string $\s(f)$ to have symbol $+$ as the $i$-th letter if $f(i-1) > 
f(i-2)$, and symbol $-$ as the $i$-th letter if $f(i-1) < f(i-2)$.  
\enddefinition

\example{Example \secea.7}
If $f$ is the path in Example \secea.4, then $\s(f) = (++--+)$.
\endexample

\proclaim{Proposition \secea.8}
The action of $W(B_l^{(1)})$ on $\Be$ induces a transitive and faithful 
action of the finite group $W_0$ of type $B_l$ on $\pml$.
For $1 \leq i < l$, the generator $s_i$ 
acts by exchanging the symbols at positions $i$ and $i+1$.  The generator $s_l$
acts by altering the symbol at position $l$.

This action extends by Theorem \secbb.5 (iv) 
to a faithful action of $W$ on $\zed \times \pml$, in which
$s_i$ acts as $(\id, s_i)$ and $s_0$ acts as follows, where $\bs'$ is the
substring obtained from $\bs$ by deleting the first two symbols: $$\eqalign{
s_0((t,\ ++ \bs')) &= (t+1,\ --\bs'),\cr
s_0((t,\ -- \bs')) &= (t-1,\ ++\bs'),\cr
s_0((t,\ +- \bs')) &= (t,\ +-\bs'),\cr
s_0((t,\ -+ \bs')) &= (t,\ -+\bs').\cr
}$$\endproclaim

\demo{Proof}
The action of $W_0$ on the strings is induced by combining the identifications
of Lemma \secbb.4 (vi) and (vii), Lemma \secea.3 and Definition 
\secea.6.  It is transitive by Lemma \secbb.4 (vi) and faithful by Theorem
\secbb.1 (iii).  The other assertions now follow from the definitions.
\qed\enddemo

\example{Example \secea.9}
Consider the string $(++--+)$ in Example \secea.7.  The action of $s_3$
on the string exchanges two minus signs, so is the identity; similarly, the
action of $s_3$ on the ideal $G$ of Example \secea.4 is the identity.
We have $s_4 . (++--+) = (++-+-)$: the symbols at positions $4$ and $5$
are exchanged.  The action of $s_4$ on $G$ adds the element $4(1)$
produce the ideal corresponding to $(++-+-)$.  We have $s_5 . (++--+)
= (++---)$: in this case the last symbol is changed.  This corresponds
to removing the element $5(1)$ from the ideal $G$.
\endexample

\head \S\seceb.  Other binary path representations \endhead

The approach of \S\secea\  can be imitated for Weyl groups of 
types $D_l^{(1)}$,
$A_{2l-1}^{(1)}$, and $C_l^{(1)}$.  Since only minor modifications to 
the arguments and definitions are required, we give only a summary of 
the results.

\topcaption{Figure \seceb.1} A full heap, $E$, over the Dynkin diagram of type 
$D_l^{(1)}$ for $l = 6$, corresponding to a spin representation,
and the convex subheap $E_0$
\endcaption
\centerline{
\hbox to 3.347in{
\vbox to 2.527in{\vfill
        \includegraphics{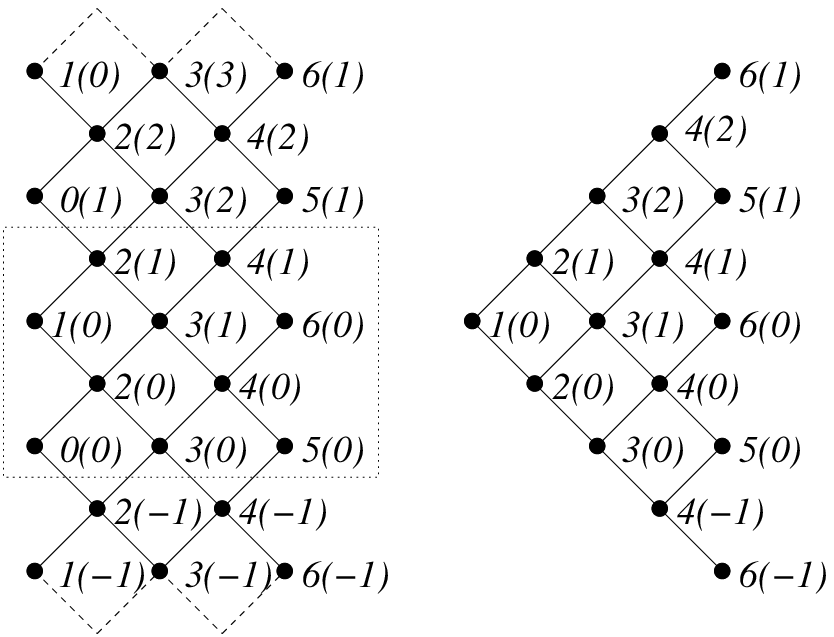}
}
\hfill}
}

The self-dual full heap $E$ in Figure \seceb.1 over the Dynkin diagram of type 
$D_l^{(1)}$ (see Figure \sectc.8) corresponds to a spin representation of the 
simple Lie algebra of type $D_l$.  The heap $E$ is ranked (see
Definition \secaa.3); we will call the 
subheap of $E$ given by the elements of rank $k$ the ``$k$-th layer'' of $E$.  
The even numbered labels in the set $X = \{2, 3, 4, \ldots, l-2\}$ occur
in the $k$-th layer if and only if $k$ is even, and the odd numbered labels
in $X$ occur in the $k$-th layer if and only if $k$ is odd.  
Label $0$ (respectively, $1$, $l$, $l-1$) occurs in the $k$-th
layer if and only if $k$ is congruent to $1$ (respectively, $3$, $l+1$,
$l+3$) modulo $4$.  (This condition ensures that the unique maximal element
of the convex subheap $E_0$ is labelled $l$.)
Another isomorphism class of heaps may be obtained in each case by twisting
by the graph automorphism exchanging vertices $l-1$ and $l$.

The ideals of the subheap $E_0$ in Figure \seceb.1 are indexed by $\pmlm$, 
similarly to the case of type $B$.
There is also an analogous notion of binary paths of type $D_l^{(1)}$.

With minor modifications to take account of the alternating pattern of
$l-1$ and $l$ in the column $x = l-2$ of the heap $E_0$, the argument 
used to prove Proposition \secea.8 proves the following

\proclaim{Lemma \seceb.1}
The action of $W(D_l^{(1)})$ on $\Be$ induces a transitive and faithful 
action of the finite group $W_0$ of type $D_l$ on $\pmlm$.
For $1 \leq i < l-1$, the generator $s_i$ acts 
by exchanging the symbols at positions $i$ and $i+1$.  Let $G$ be the ideal
of $E_0$ corresponding to a string $\bs$ of length $l-1$, and define
$k = g(G, l-2)$, where $g$ is as in Definition \secea.1.  Call $G$ ``even''
if $k = l-1 \mod 4$, and ``odd'' otherwise (\idest if $k = l+1 \mod 4$).
Then $s_{l-1}$ and $s_l$ act on the $(l-1)$-st symbol of $\bs$ according to
the following rules: $$\eqalign{
s_{l-1}(-) &= \cases + & \text{ if } G \text{ is even},\cr
- & \text{ otherwise;}\cr
\endcases\cr
s_{l-1}(+) &= \cases - & \text{ if } G \text{ is odd},\cr
+ & \text{ otherwise,}\cr
\endcases\cr
s_{l}(-) &= \cases + & \text{ if } G \text{ is odd},\cr
- & \text{ otherwise;}\cr
\endcases\cr
s_{l}(+) &= \cases - & \text{ if } G \text{ is even},\cr
+ & \text{ otherwise.} \qed\cr
\endcases} $$  
\endproclaim

Lemma \seceb.1 can be summarised more concisely by appending a symbol ``$+$''
to $\bs$ if $\bs$ corresponds to an even ideal, and a symbol ``$-$'' if it
corresponds to an odd ideal.  Note that this will produce elements of
$\pml$ containing an even number of $-$ signs.

\proclaim{Proposition \seceb.2}
The action of $W(D_l^{(1)})$ on $\Be$ induces a transitive and faithful 
action of the finite group $W_0$ of type $D_l$ on the subset of $\pml$
consisting of strings that contain an even number of $-$ signs.  
For $1 \leq i < l$, the generator $s_i$ 
acts by exchanging the symbols at positions $i$ and $i+1$.  The generator $s_l$
acts by exchanging the symbols at positions $l-1$ and $l$, and then altering
each of the symbols at these positions.
This action extends by Theorem \secbb.5 (iv) 
to a faithful action of $W$ on $\zed \times \pml$,
in which $s_i$ acts as $(\id, s_i)$ for $i \ne 0$, 
and the action of $s_0$ is as described in Proposition \secea.8.
\qed\endproclaim

We now turn to type $A_{2l-1}^{(1)}$, whose Dynkin diagram is shown in
Figure \sectc.1; note that the diagram has an even number, $2l$, of vertices.
The relevant full heap is the ranked heap $E$ for which the
$k$-th layer consists of all odd (respectively, even) numbered vertices 
if $k$ is odd (respectively, even).  It is convenient to regard the
elements of $E$ as occuping places in a $\zed_{2l} \times \zed$ grid on a
cylinder, in which vertex $p$ in layer $k$ corresponds to the position
$(p, k)$.  (Note that precisely half the possible positions are occupied
by elements of $E$.) 

\definition{Definition \seceb.3}
To each proper ideal $F$ of the heap $E$ defined above, we define $$
g(F, i) = \max \{y : (i, y) \in F\} + 1
.$$  We also define
a binary path of type $A_{2l-1}^{(1)}$ to be a function $$
f : \zed_{2l} \ra \zed
$$ such that for all $i \in \zed$ we have 
$f(i+1) = f(i) \pm 1$, where indices are read modulo $2l$.  We denote the 
set of all binary paths of type $A_{2l-1}^{(1)}$ by $P(A_{2l-1}^{(1)})$.
\enddefinition

We omit the proof of the following lemma, because it is very similar to 
the proof of Lemma \secea.3.

\proclaim{Lemma \seceb.4}
There is a bijection $\pi$ from the set of binary paths of type 
$A_{2l-1}^{(1)}$
to the set of proper ideals of the heap $E$ defined above,
identifying a proper ideal with the corresponding subset of $\zed_{2l} 
\times \zed$ as above.  Explicitly, we have $$
\pi(f) = \{(x, y) \in E : y < f(x)\}
.\qed$$\endproclaim

The corresponding notion of strings in this case is the set $\pmc$, which
consists of those functions $f : \zed_{2l} \ra \{+, -\}$ such that 
$|f^{-1}(+)| = |f^{-1}(-)|$.

\definition{Definition \seceb.5}
Given a binary path $f$, we define
the string $\s(f) \in \pmc$ to have symbol $+$ as the $i$-th letter 
if $f(i) > f(i-1)$, and symbol $-$ as the $i$-th letter if $f(i) < f(i-1)$,
where indices are read modulo $2l$.
(Note that $\s(f)$ will have an equal number of $+$ and $-$ signs.)
\enddefinition

The analogue of propositions \secea.8 and \seceb.2 for type $A_{2l-1}^{(1)}$
is as follows; the proof follows the same lines as Proposition \secea.8.
In order to state the result, it is convenient to represent an element $\bs
\in \pmc$ by the sequence $(f(0), f(1), \ldots, f(2l-1))$, and to define 
$\bs'$ to be the substring obtained from $\bs$ by deleting the first and last 
elements in the sequence representation.  

\proclaim{Proposition \seceb.6}
The action of $W(A_{2l-1}^{(1)})$ on $\Be$ induces a transitive and faithful 
action of the finite group $W_0$ of type $A_{2l-1}$ on $\pmc$.
For $i \not\equiv 0 \mod 2l$, the generator $s_i$ 
acts by exchanging the symbols at positions $i-1$ and $i$.
This action extends by Theorem \secbb.5 (iv) 
to a faithful action of $W$ on $\zed \times \pmc$,
in which $s_i$ acts as $(\id, s_i)$ for $i \not\equiv 0 \mod 2l$, 
and $s_0$ acts as follows: $$\eqalign{
s_0((t,\ + \bs' -)) &= (t+1,\ - \bs' +),\cr
s_0((t,\ - \bs' +)) &= (t-1,\ + \bs' -),\cr
s_0((t,\ + \bs' +)) &= (t,\ + \bs' +),\cr
s_0((t,\ - \bs' -)) &= (t,\ - \bs' -). \qed\cr
}$$\endproclaim

\topcaption{Figure \seceb.2} The Dynkin diagram of type $D_{l+1}^{(2)}$
\endcaption
\centerline{
\hbox to 3.388in{
\vbox to 0.388in{\vfill
        \includegraphics{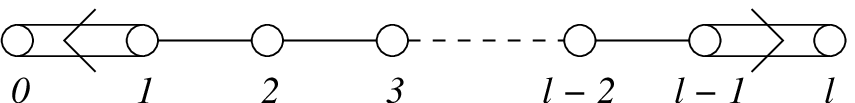}
}
\hfill}
}

There is also a binary path representation for the affine Weyl group 
$W(C_l^{(1)}) \cong W(D_{l+1}^{(2)})$ of
type $C$.
The corresponding full
heap is similar to the one for type $B_l^{(1)}$ (see Figure \secea.2), except
that the generators corresponding to $0$ and $1$ in Figure \secea.1 are 
identified with each other.  More precisely, the relevant full heap 
is the ranked heap $E$ over the Dynkin diagram of Figure \seceb.2 for which the
$k$-th layer consists of all odd (respectively, even) numbered vertices 
if $k$ is odd (respectively, even).  
The Weyl group action may be explicitly described as follows.

\proclaim{Proposition \seceb.7}
The action of $W(C_l^{(1)})$ on $\Be$ induces a transitive and faithful 
action of the finite group $W_0$ of type $C_l$ on $\pmlp$.
For $i \ne 0$, the generator $s_i$ 
acts by exchanging the symbols at positions $i$ and $i+1$.
This action extends by Theorem \secbb.5 (iv) 
to a faithful action of $W$ on $\zed \times \pmlp$,
in which $s_i$ acts as $(\id, s_i)$ for $i \ne 0$, 
and $s_0$ acts as follows, where $\bs'$ is the substring obtained from $\bs$
by deleting the first symbol: $$\eqalign{
s_0((t,\ + \bs')) &= (t+1,\ - \bs'),\cr
s_0((t,\ - \bs')) &= (t-1,\ + \bs').\qed\cr
}$$
\endproclaim

\remark{Remark \seceb.8}
The representations described in propositions \secea.8, \seceb.2, \seceb.7
and \seceb.6
may, by taking a suitable limit, be extended to representations of the
infinitely generated Weyl groups of types $A_\infty$, $B_\infty$, 
$C_\infty$ and $D_\infty$ respectively.  (The Weyl groups of types $B_\infty$
and $C_\infty$ are isomorphic.  See \cite{{\bf 10}, \S7.11} for the relevant
definitions.)
In these cases, the full heaps over the Dynkin diagrams are countably
infinite as before, but the lattice $\Be$ is uncountable.
It can be shown that the Weyl group acts faithfully in these cases,
although of course the action is intransitive, for reasons of cardinality.  We 
omit the details for reasons of space.
\endremark

\head \S\sectf. Some related constructions \endhead

In \S\sectf, we look briefly at some other examples of permutation 
representations of Weyl groups and related groups, in particular,
examples arising from 
the lattices of ideals of heaps that are not full heaps.  We will not give
full details in the interests of space.

A well known result in the theory of lattices \cite{{\bf 16}, \S3.4} is that
the set of ideals of a partially ordered set $(E, \leq)$ forms 
a distributive lattice, $J = J(E)$, where the operations of meet and join are
set theoretic intersection and union, respectively.  A covering relation
in $J$ consists of two ideals $F \subset F'$ such that $F' \backslash F$
is a singleton.  In the case where
$(E, \leq)$ is a heap, the elements of $E$ have labels, so in the case
of a covering relation, we may label the edge from $F$ to $F'$ in the Hasse
diagram for $J$ by the label of the single element in $F' \backslash F$.

In \cite{{\bf 7}, Theorem 8.3}, it is shown how the distributive lattice associated
to a full heap over an (untwisted) affine Kac--Moody algebra has the structure
of a crystal in the sense of Kashiwara \cite{{\bf 11}}; both these structures
can be regarded as edge-labelled Hasse diagrams of partially ordered sets.

Another class of examples of edge-labelled graphs in the theory of Coxeter 
groups are the so-called Cox Box blocks of Eriksson's thesis \cite{{\bf 6}, \S2.7}.
Cox Box blocks are easy to construct, although there is a lot of freedom of
choice in their construction; as Eriksson says in \cite{{\bf 6}, \S2.7}, they are
``perhaps more of a game than a systematic approach''.   
Another disadvantage of Cox Box blocks is that they
need not lead to faithful representations of the associated Coxeter group
\cite{{\bf 6}, \S2.7.1}.
Remarkably however, although the Cox Box blocks have no inherent partial 
order associated to them, many of the most natural examples do turn out to be
isomorphic as edge-labelled graphs to the lattice of ideals of a heap.
Moreover, in many cases, the heap that arises is a full heap, which by 
the theory in this paper
shows that the associated Coxeter group representation will be faithful, at
least if 
the associated Kac--Moody algebra is one of the affine types.
From this point of view, the action of a Coxeter generator, $s_i$, on an ideal
is easy to describe directly, as follows.  First locate the vertex $v$ of
the edge-labelled graph corresponding to the ideal.  If there is no edge
labelled $i$ emerging from $v$, then $s_i$ acts as the identity.  Otherwise,
follow the (unique) edge labelled $i$ to another vertex, $v'$ say; then
$s_i . v = v'$.

Some examples of Cox Box blocks are shown in \cite{{\bf 6}, p62}.
The type $F_4$ graph given there may not be the Hasse diagram of a 
distributive lattice, but the other ten examples shown all occur as
the lattice of ideals of some heap, and in all but two cases, this heap is a
full heap.  These two exceptions are types $G_2^{(1)}$ and type $H_3$,
as we now explain.

\vfill\eject

\topcaption{Figure \sectf.1} A full heap of type $G_2^{(1)}$ (left), and its
corresponding lattice of ideals (right)
\endcaption
\centerline{
\hbox to 2.680in{
\vbox to 2.500in{\vfill
        \includegraphics{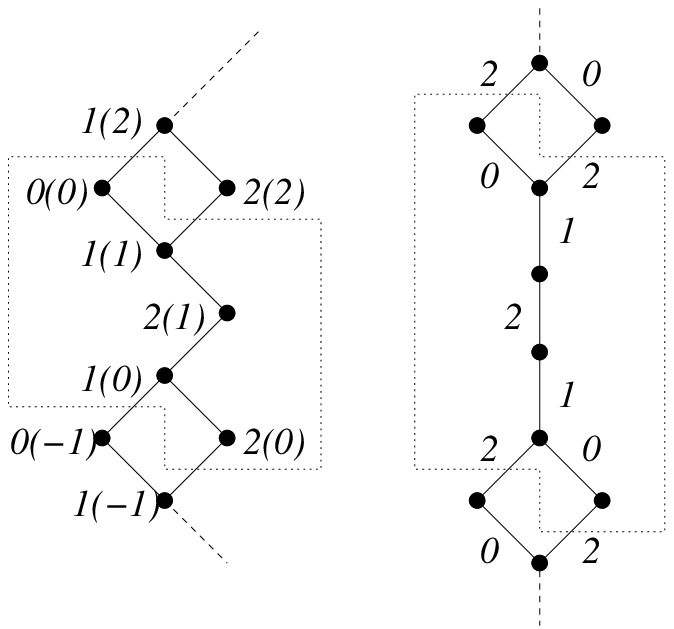}
}
\hfill}
}

\example{Example \sectf.1}
Let $W$ be the affine Weyl group of type $G_2^{(1)}$ given by the 
presentation $$
W = \lan s_0, s_1, s_2 : s_0^2 = s_1^2 = s_2^2 = 1, 
(s_0 s_1)^3 = 1, (s_1 s_2)^6 = 1, (s_0 s_2)^2 = 1
.\ran
$$  Ignoring the labels on the edges, the Dynkin diagram in this case is
a graph $\Gamma$ in which $s_1$ is connected to $s_0$ and $s_2$, and $s_0$ 
is not connected to $s_2$.  Consider the heap $E$ over $\Gamma$ shown on the
left in Figure \sectf.1.  The Hasse diagram of the lattice of ideals of $E$, 
considered as a labelled graph, is shown on the right in Figure \sectf.1.  This
lattice is isomorphic as a labelled graph to Eriksson's Cox Box blocks
model for type $G_2^{(1)}$ \cite{{\bf 6}, p62}.

Remarkably, our techniques can be modified to produce a categorification
of the action of the Weyl group on the root system in this case.  The key
to this is to think of the heap elements $E(2, 2k+1)$ as ``double vertices'':
more precisely, we redefine $\cha$ so that $\cha(\{E(2, 2k+1)\}) = 2\a_2$,
but we let $\cha(\{E(2, 2k)\}) = \a_2$ as before.  Our arguments can then
be adapted to prove that the representation is faithful.  This is
not isomorphic to the representation constructed by Cellini et al in 
\cite{{\bf 4}, Theorem 4.5}, as the latter contains points that are fixed by 
all Weyl group elements, although it may be the case that the representations
become isomorphic after the removal of these fixed points.
It may also be possible to construct
a Lie algebra representation from our example by ad hoc modification of 
the definitions, as Wildberger does in \cite{{\bf 20}} for finite type
$G_2$.
\endexample

\topcaption{Figure \sectf.2} A full heap of type $H_3$ (left), and its
corresponding lattice of ideals (right)
\endcaption
\centerline{
\hbox to 2.388in{
\vbox to 2.583in{\vfill
        \includegraphics{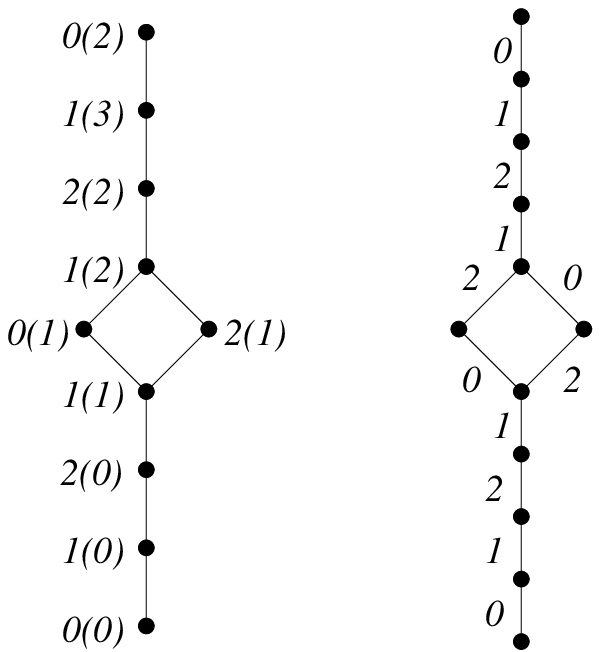}
}
\hfill}
}

\example{Example \sectf.2}
Let $W$ be the finite Coxeter group of type $H_3$ given by the 
presentation $$
W = \lan s_0, s_1, s_2 : s_0^2 = s_1^2 = s_2^2 = 1, 
(s_0 s_1)^3 = 1, (s_1 s_2)^5 = 1, (s_0 s_2)^2 = 1
.\ran
$$  Ignoring the labels on the edges, the Dynkin diagram in this case is
a graph $\Gamma$ in which $s_1$ is connected to $s_0$ and $s_2$, and $s_0$ 
is not connected to $s_2$.  Consider the heap $E$ over $\Gamma$ shown on the
left in Figure \sectf.2.  The Hasse diagram of the lattice of ideals of $E$, 
considered as a labelled graph, is shown on the right in Figure \sectf.2.  This
lattice is isomorphic as a labelled graph to Eriksson's Cox Box blocks
model for type $H_3$ \cite{{\bf 6}, p62}.

Eriksson \cite{{\bf 6}, \S2.7.2} states (without proof, but correctly!) 
that this representation
is permutation group isomorphic to the action of the Coxeter group of type
$H_3$ as rigid rotations and reflections on the twelve vertices of
the icosahedron; it follows that the representation is faithful as both
the Coxeter group and the symmetry group are known to have order $120$.
\endexample

\remark{Remark \sectf.3}
Note that the edge-labelled graphs in examples \sectf.1 and \sectf.2 are not 
much more complicated than the heaps that give rise to them: this is
because the heap in each case 
is close to being totally ordered.  In other examples,
such as those of \S\secea\  and \S\seceb, the Cox Box blocks construction is 
so much more
complicated than the heap that it becomes difficult to write down; an extreme
example of this is the situation of Remark \seceb.8, in which the 
edge-labelled graph has uncountably many vertices.
\endremark

\head Acknowledgement \endhead

I thank H. Eriksson for sending me a copy of \cite{{\bf 6}}.

\leftheadtext{} \rightheadtext{}

\end